\documentclass[11pt]{amsart}
\usepackage{amssymb,latexsym,amsmath,graphicx,graphics,epic,eepic}

\newcommand{\zed}{\mathbb{Z}}

\theoremstyle{plain}
\newtheorem{theorem}{Theorem}
\numberwithin{theorem}{section}
\newtheorem{lemma}[theorem]{Lemma}
\newtheorem{proposition}[theorem]{Proposition}
\newtheorem{corollary}[theorem]{Corollary}
\newtheorem{conjecture}[theorem]{Conjecture}

\theoremstyle{definition}
\newtheorem{definition}[theorem]{Definition}

\newtheorem{remark}[theorem]{Remark}
\newtheorem{acknowledgments}{Acknowledgments\ignorespaces}

\newcommand{\id}{\text{id}}


\begin{document}

\title{Braids, Transversal links and the Khovanov-Rozansky Theory}

\author{Hao Wu}

\address{Department of mathematics and Statistics\\
         Lederle Graduate Research Tower\\
         710 North Pleasant Street\\
         University of Massachusetts\\
         Amherst, MA 01003-9305\\
         USA}

\email{wu@math.umass.edu}

\begin{abstract}
We establish some inequalities about the Khovanov-Rozansky
cohomologies of braids. These give new upper bounds of the
self-linking numbers of transversal links in standard contact
$S^3$ which is sharper than the well known bound given by the
HOMFLY polynomial. We also introduce a sequence of transversal
link invariants, and discuss some of their properties.
\end{abstract}

\maketitle

\section{Introduction}

In \cite{KR1,KR2},  M. Khovanov and L. Rozansky introduced two
versions of Khovanov-Rozansky link cohomologies. Their construction is
based on matrix factorizations associated to certain planar
diagrams. The version in \cite{KR1} uses the potential
$x^{n+1}$, which lead to a $\mathbb{Z}\oplus\mathbb{Z}$-graded
cohomology theory $H_n$ for each $n\geq2$. The version in \cite{KR2}
uses the potential $ax$, which gives a
$\mathbb{Z}\oplus\mathbb{Z}\oplus\mathbb{Z}$-graded cohomology
theory $H$.

M. Khovanov and L. Rozansky showed that $H_n$ is invariant under
Reidemeister moves, and the
isomorphism type of $H$ is invariant under Markov moves of closed
braids up to overall shifts of gradings. The graded Euler
characteristics of $H_n$ and $H$ are variants of the HOMFLY
polynomial.

\begin{figure}[h]

\setlength{\unitlength}{1pt}

\begin{picture}(420,50)(-210,-30)

\linethickness{.5pt}


\put(-100,-20){\vector(1,1){40}}

\put(-60,-20){\line(-1,1){15}}

\put(-85,5){\vector(-1,1){15}}

\put(-84,-30){$+$}

\put(-115,0){$\sigma_i$}


\put(100,-20){\vector(-1,1){40}}

\put(60,-20){\line(1,1){15}}

\put(85,5){\vector(1,1){15}}

\put(76,-30){$-$}

\put(45,0){$\sigma_i^{-1}$}

\end{picture}

\caption{Sign of a crossing}\label{sign}

\end{figure}
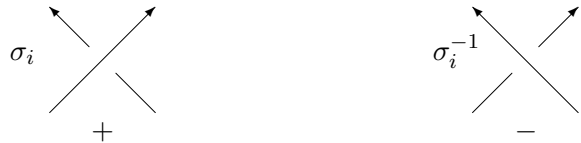

In this paper, we discuss some applications of these theories in
contact topology, specifically the transversal knot theory. The
sign convention of this paper is given in figure \ref{sign}, which
is opposite to that used in \cite{KR1,KR2}.

First, we introduce a variant of $H$ whose graded isomorphism type
is invariant under Reidemeister moves. Let $B$ be a closed braid
with $b$ strands and writhe $w$. Define $\hat{H}(B)$ to be $H(B)$
with each of the two quantum gradings shifted by $\frac{w+b}{2}$ and the homological grading shifted by $-\frac{w+b}{2}$. Then
$\hat{H}(B)$ is
$\frac{1}{2}\mathbb{Z}\oplus\frac{1}{2}\mathbb{Z}\oplus\frac{1}{2}\mathbb{Z}$-graded.
From M. Khovanov and L. Rozansky's construction, it's clear that the
gradings of $\hat{H}(B)$ are now invariant under all Markov moves.

Denote by $H_n^{j,k}$ the space of homogeneous elements with
grading $(j,k)$ of $H_n$, where $j$ is the homological grading,
and $k$ is the quantum grading which comes from the exponents of
the markings $x_i$. Denote by $\hat{H}^{j,k,l}(B)$ the space of
homogeneous elements with grading $(j,k,l)$ of $\hat{H}(B)$, where
$j$ is the cohomological grading, $k$ is the quantum $a$-grading
which comes from the exponents of the formal variable $a$, and $l$
is the quantum $x$-grading which comes from the exponents of
markings $x_i$.

For a link $L$, define
\[
g^{(n)}_{max}(L)=\max\{k~|~H_n^{j,k}(L)\neq0 \text{ for some }j\},
\]
\[
g^{(n)}_{min}(L)=\min\{k~|~H_n^{j,k}(L)\neq0 \text{ for some }j\},
\]
\[
\hat{g}_{max}(L)=\max\{k~|~\hat{H}^{j,k,l}(L)\neq0 \text{ for some
}j,~l\},
\]
\[
\hat{g}_{min}(L)=\min\{k~|~\hat{H}^{j,k,l}(L)\neq0 \text{ for some
}j,~l\}.
\]
These are numerical link invariants. We have:

\begin{theorem}\label{main}
If $B$ is a closed braid with $b$ strands and writhe $w$, then
\begin{equation}\label{main-inequality}
\frac{w-b}{2} \leq \hat{g}_{min}(B) \leq \hat{g}_{max}(B) \leq
\frac{w+b}{2}-1.
\end{equation}
If we further assume that the numbers of positive and negative
crossings of $B$ are $c_+$ and $c_-$, respectively, then
\begin{equation}\label{main-inequality-n}
(n-1)(w-b)-2c_- \leq g^{(n)}_{min}(B) \leq g^{(n)}_{max}(B) \leq
(n-1)(w+b)+2c_+,
\end{equation}
\end{theorem}

\begin{definition}
The braid number $b(L)$ of a link $L$ is the minimal number of
strands necessary to represent $L$ as a braid.
\end{definition}

\begin{corollary}\label{braid-number}
For a link $L$,
\[
b(L)\geq \hat{g}_{max}(L)- \hat{g}_{min}(L) +1,
\]
\[
b(L)\geq \frac{1}{2}\limsup_{n\rightarrow\infty}
\frac{g^{(n)}_{max}(L)-g^{(n)}_{min}(L)}{n-1}.
\]
\end{corollary}
\noindent\textit{Proof.} $\Box$

\begin{theorem}\label{tb}
If $L$ is a transversal link in the standard contact $3$-sphere,
then
\[
sl(L)\leq 2\hat{g}_{min}(L),
\]
\[
sl(L) \leq -\limsup_{n\rightarrow\infty}
\frac{g^{(n)}_{max}(\overline{L})}{n-1},
\]
where "$sl$" is the self-linking number, and $\overline{L}$ is the
mirror image of $L$.

Therefore, if $K$ is a Legendrian knot in the standard contact
$3$-sphere, then
\[
tb(K)+|r(K)|\leq 2\hat{g}_{min}(K),
\]
\[
tb(K)+|r(K)| \leq -\limsup_{n\rightarrow\infty}
\frac{g^{(n)}_{max}(\overline{K})}{n-1},
\]
where "$tb$" is the Thurston-Bennequin number, and "$r$" is the rotation
number.
\end{theorem}
\begin{proof}
By \cite{Ben}, any transversal knot is transversally isotopic to a
transversal braid, and the self-linking number $sl(B)$ of a
transversal braid $B$ with $b$ strands is equal to $w-b$, where
$w$ is the writhe of $B$. So Theorem \ref{main} implies the first
part of the theorem. The second part follows from the first part
by a push-off argument given in \cite{Ben}.
\end{proof}

In \cite{Pl4}, O. Plamenevskaya constructed an invariant $\psi$ of
transversal links in standard contact $S^3$ using the original
Khovanov homology $\mathcal{H}$ defined in \cite{K1}. We
generalize her construction, and, for each $n\geq2$, define an
invariant $\psi_n$ of transversal links in standard contact $S^3$
using $H_n$.

\begin{theorem}\label{invariants}
Let $L$ be a transversal link in the standard contact $S^3$. For
each $n\geq2$, we can associate to $L$ an element $\psi_n(L)$ of
$H_n(\overline{L})$ of cohomological degree $0$ and quantum
degree $-(n-1)sl(L)$. $\psi_n(L)$ is invariant under transversal
isotopy of $L$ up to multiplication by a non-zero scalar.

When $n=2$, $\psi_2(L)$ is identified with $\psi(L)$ under the
isomorphism
\[
H_2(\overline{L})\cong
\mathcal{H}(L)\otimes_{\zed}\mathbb{Q}.
\]
\end{theorem}

Comparing Theorems \ref{tb} and \ref{invariants}, it's easy to see
that the sequence $\{\psi_n\}$ sometimes detects the maximal
self-linking number. Specially, we can generalize Plamenevskaya's
computation about quasi-positive braids, and get a stronger
conclusion.

\begin{corollary}\label{top-sl}
If $L$ is a transversal link in the standard contact $S^3$ with
$\psi_n(L)\neq0$ for infinitely many $n$, then
\[
SL(L)=sl(L)=-\limsup_{n\rightarrow\infty}
\frac{g^{(n)}_{max}(\overline{L})}{n-1},
\]
where $SL(L)$ is the maximal self-linking number of
all transversal links smoothly isotopic to $L$.

Specially, if a transversal link $L$ is transverse isotopic to a
quasi-positive transversal braid, then $\psi_n(L)\neq0$, $\forall
~n\geq2$, and, thus, the above equation is true.
\end{corollary}

\begin{proof}
The first part of the corollary is clear from Theorems \ref{tb}
and \ref{invariants}. The second part will be proved in
Proposition \ref{quasi-positive}.
\end{proof}

A transversal knot invariant is said to be non-classical if it is
not determined by the smooth knot type and the self-linking
number. There are still no known examples of efficient
non-classical transversal knot invariants.

\begin{conjecture}
The sequence $\{\psi_n\}$ is a non-classical transversal knot
invariant.
\end{conjecture}

The rest of this paper is organized as following: In Section
\ref{kr}, we briefly review the definitions of the
Khovanov-Rozansky cohomologies and compare the graded Euler
characteristics of these cohomologies with the HOMFLY polynomial,
which explains why Theorems \ref{main}, \ref{tb} and Corollary
\ref{braid-number} imply the corresponding results by J. Franks,
R. Williams and H. Morton in \cite{FW,Mo}, including the upper
bound of the self-linking number from the HOMFLY polynomial. In
Section \ref{inequality}, we study the Khovanov-Rozansky
cohomologies of closed braids, and prove Theorem \ref{main}. In
Section \ref{inv}, we construct the invariants $\psi_n$, and study
some of their properties.

\section{Khovanov-Rozansky Theory}\label{kr}

\subsection{Matrix factorizations}
Let $R=\mathbb{Q}[t_1,\cdots,t_m]$ be a polynomial ring.  A matrix
factorization $M$ over $R$ with potential $w\in R$ is a collection
of two $R$-modules $M^0$, $M^1$ and two module homomorphisms
$d^0:M^0\rightarrow M^1$, $d^1:M^1\rightarrow M^0$, called
differential maps, s.t.,
\[
d^1d^0=w\cdot\id_{M^0}, \hspace{1cm} d^0d^1=w\cdot\id_{M^1}.
\]
We usually write such a matrix factorization $M$ as
\[
M^0 \xrightarrow{d_0} M^1 \xrightarrow{d_1} M^0.
\]

For $a,b\in R$, denote by $(a,b)_R$ the matrix factorization
\[
R \xrightarrow{a} R \xrightarrow{b} R,
\]
where $a,~b$ act on $R$ by multiplication. This matrix
factorization has potential $ab$. When the polynomial ring $R$ is
clear from context, we drop it from the notation, and only write
$(a,b)$.

Denote by
\[
\left(%
\begin{array}{cc}
  a_1 & b_1 \\
  a_2 & b_2 \\
  \vdots & \vdots \\
  a_k & b_k \\
\end{array}%
\right)_R
\]
the tensor product of
$(a_1,b_1)_R,~(a_2,b_2)_R,\cdots,~(a_k,b_k)_R$, where the
differential maps are given by the (signed) Leibniz rule. Its
potential is
\[
w=a_1b_1+a_2b_2+\cdots+a_kb_k.
\]
Again, when $R$ is clear from context, we drop it from the
notation.

\subsection{Grading shifts}
Suppose $M$ is a $\mathbb{Z}$-graded module over a
$\mathbb{Z}$-graded ring, and $m\in\mathbb{Z}$. Define $M\{m\}$ to
be $M$ with grading shifted by $m$, i.e., an element $x$ of
$M\{m\}$ is homogeneous of degree $i+m$ if and only if it's a
homogeneous element of $M$ of degree $i$. Similarly, suppose $M$
is a $\mathbb{Z}\oplus\mathbb{Z}$-graded module over a
$\mathbb{Z}\oplus\mathbb{Z}$-graded ring, and $m,l\in\mathbb{Z}$.
Define $M\{m,l\}$ to be $M$ with grading shifted by $(m,l)$, i.e.,
an element $x$ of $M\{m,l\}$ is homogeneous of bi-degree $(i+m,
j+l)$ if and only if it's a homogeneous element of $M$ of degree
$(i,j)$.

\subsection{Planar diagrams}\label{pd}
To construct the Khovanov-Rozansky cohomology theories, one
considers the planar diagrams $\Gamma$ with the following
properties:
\begin{enumerate}
    \item $\Gamma$ consists of two types of edges: regular edges
    and wide edges. These edges intersect only at their endpoints.
    \item Regular edges are disjoint from each other.
    Wide edges are disjoint from each other.
    \item Each regular edge is oriented, and contains at least one
    marked point. Open endpoints of regular edges are marked.
    \item Each wide edge has exactly two regular edges entering at one endpoint,
    and exactly two regular edges exiting from the other endpoint.
\end{enumerate}

See Figures \ref{eg}, \ref{reduction1}, \ref{reduction2}, \ref{reduction3} below for examples
of such diagrams.

\subsection{Matrix factorization associated to a planar diagram}
Let $\Gamma$ be a planar diagram satisfying the conditions in
subsection \ref{pd}, and $\{x_1, \cdots, x_p\}$ the set of
markings on $\Gamma$. Let $R'=\mathbb{Q}[x_1,\cdots,x_p]$ be the
$\mathbb{Z}$-graded polynomial ring so that $\deg x_i=2$,
$\forall~i$. Let $R=\mathbb{Q}[a,x_1,\cdots,x_p]$ be the
$\mathbb{Z}\oplus\mathbb{Z}$-graded polynomial ring so that $\deg
a=(2,0)$ and $\deg x_i=(0,2)$, $\forall~i$.

\begin{figure}[h]

\setlength{\unitlength}{1pt}

\begin{picture}(420,20)(-210,0)

\linethickness{.5pt}

\put(-15,0){\vector(1,0){30}}

\put(-20,5){$x_j$}

\put(15,5){$x_i$}

\end{picture}

\caption{$L^i_j$}\label{lij}

\end{figure}

For an oriented (regular) arc $L^i_j$ in $\Gamma$ from the point
marked by $x_j$ to the point marked by $x_i$ that has no marked
interior points, let $C_n(L^i_j)$ be the matrix factorization
$(\pi_{ij},x_i-x_j)_{R'}$ given by
\[
R' \xrightarrow{\pi_{ij}} R'\{1-n\} \xrightarrow{x_i-x_j} R',
\]
where $\pi_{ij}=x_i^n+x_i^{n-1}x_j+\cdots+x_ix_j^{n-1}+x_j^n$. The
purpose of the grading shift here is to make $C_n(L^i_j)$ graded
in the sense that both $\pi_{ij}$ and $x_i-x_j$ are homogeneous
maps of degree $n+1$.

Also, let $C(L^i_j)$ be the matrix factorization $(a,x_i-x_j)_R$
given by
\[
R \xrightarrow{a} R\{-1,1\} \xrightarrow{x_i-x_j} R,
\]
where the grading shift makes it bi-graded in the sense that $a$
and $x_i-x_j$ are both homogeneous maps of bi-degree $(1,1)$.

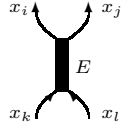
\begin{figure}[h]

\setlength{\unitlength}{1pt}

\begin{picture}(420,45)(-210,0)

\linethickness{.5pt}

\put(-2.5,10){\vector(1,1){0}}

\qbezier(-10,0)(-10,5)(-2.5,10)

\qbezier(-10,40)(-10,35)(-2.5,30)

\put(-10,44){\vector(0,1){0}}

\put(2.5,10){\vector(-1,1){0}}

\qbezier(10,0)(10,5)(2.5,10)

\qbezier(10,40)(10,35)(2.5,30)

\put(10,44){\vector(0,1){0}}

\put(-20,40){\tiny{$x_i$}}

\put(15,40){\tiny{$x_j$}}

\put(-20,0){\tiny{$x_k$}}

\put(15,0){\tiny{$x_l$}}

\linethickness{5pt}

\put(0,10){\vector(0,1){20}}

\put(5,17){\tiny{$E$}}

\end{picture}

\caption{Wide edge $E$}\label{wideedge}

\end{figure}

For a wide edge $E$, let $x_i$, $x_j$, $x_k$, $x_l$ be the closest
markings to $E$ as depicted in Figure \ref{wideedge}. (It is
possible that $x_i=x_k$ or $x_j=x_l$.) Let $g$ be the unique
two-variable polynomial such that $g(x+y,xy)=x^{n+1}+y^{n+1}$, and
\[
u_1 = u_1(x_i,x_j,x_k,x_l) =
\frac{g(x_i+x_j,x_ix_j)-g(x_k+x_l,x_ix_j)}{x_i+x_j-x_k-x_l},
\]
\[
u_1 = u_1(x_i,x_j,x_k,x_l) =
\frac{g(x_k+x_l,x_ix_j)-g(x_k+x_l,x_kx_l)}{x_ix_j-x_kx_l}.
\]
Note that $u_1$ and $u_2$ are homogeneous polynomials in $x_i$,
$x_j$, $x_k$ and $x_l$ of degrees $2n$ and $2n-2$, respectively.

Define $C_n(E)$ to be the matrix factorization
\[
\left(%
\begin{array}{cc}
  u_1 & x_i+x_j-x_k-x_l \\
  u_2 & x_ix_j-x_kx_l \\
\end{array}%
\right)_{R'}\{-1\},
\]
which is the tensor product of the matrix factorizations
\[
R'\{-1\} \xrightarrow{u_1} R'\{-n\} \xrightarrow{x_i+x_j-x_k-x_l}
R'\{-1\}
\]
and
\[
R' \xrightarrow{u_2} R'\{3-n\} \xrightarrow{x_ix_j-x_kx_l} R'.
\]
Again, the grading shifts here make $C_n(E)$ graded in the sense
that both differential maps are homogeneous maps of degree $n+1$.

Also define $C(E)$ to be the matrix factorization
\[
\left(%
\begin{array}{cc}
  a & x_i+x_j-x_k-x_l \\
  0 & x_ix_j-x_kx_l \\
\end{array}%
\right)_R,
\]
which is the tensor product of the matrix factorizations
\[
R \xrightarrow{a} R\{-1,1\} \xrightarrow{x_i+x_j-x_k-x_l} R
\]
and
\[
R \xrightarrow{0} R\{-1,3\} \xrightarrow{x_ix_j-x_kx_l} R.
\]
Again, the grading shifts here make $C(E)$ bi-graded in the sense
that both differential maps are homogeneous maps of degree
$(1,1)$.

Finally, we define that
\[
C_n(\Gamma) = (\bigotimes_{L^i_j}C_n(L^i_j)) ~\bigotimes~
(\bigotimes_{E}C_n(E)),
\]
and
\[
C(\Gamma) = (\bigotimes_{L^i_j}C(L^i_j)) ~\bigotimes~
(\bigotimes_{E}C(E)).
\]
where $L^i_j$ runs through all the (regular) oriented arcs
starting and ending at marked points with no marked interior
points, and $E$ runs through all wide edges.

If $\Gamma$ is closed, i.e., has no open end points, then
$C_n(\Gamma)$ and $C(\Gamma)$ are graded cyclic chain complexes.
We denote by $H_n(\Gamma)$ and $H(\Gamma)$ their cohomologies.
These cohomologies inherit the gradings of $C_n(\Gamma)$ and
$C(\Gamma)$. Note that the $\zed_2$-grading from the cyclic
complex structure is trivial since all the non-vanish cohomology
concentrate on one of the two degrees. The cohomologies
$H_n(\Gamma)$ and $H(\Gamma)$ do not depend on the choice of
marked points.

If $\Gamma$ has open end points, let $x_{i_1},\cdots,~x_{i_q}$ be
the markings of these open end points. Then $C_n(\Gamma)$ is a
matrix factorization of potential $\sum_{r=1}^{q} \pm x^{n+1}_{i_r}$,
and $C(\Gamma)$ is a factorization of potential $\sum_{r=1}^{q} \pm
ax_{i_r}$, where the sign is "$+$" if the corresponding open end
point is an exit, and is "$-$" if otherwise. Let $I'$ and $I$ be
the ideals generated by $\{x_{i_1},\cdots,~x_{i_q}\}$ in $R'$ and
$R$, respectively. Then $C_n(\Gamma)/I'C_n(\Gamma)$ and
$C(\Gamma)/IC(\Gamma)$ are graded cyclic chain complex. Again, we
denote by $H_n(\Gamma)$ and $H(\Gamma)$ their cohomologies, which
have the gradings inherited from $C_n(\Gamma)$ and $C(\Gamma)$. As
before, the cohomologies $H_n(\Gamma)$ and $H(\Gamma)$ do not
depend on the choice of marked points.

\subsection{Khovanov-Rozansky cohomologies of a knot}

\begin{figure}[h]

\setlength{\unitlength}{1pt}

\begin{picture}(420,60)(-210,-15)


\put(-40,0){\vector(0,1){40}}

\put(-60,0){\vector(0,1){40}}

\put(-52,-15){\small{$\Gamma_0$}}

\put(-70,40){\tiny{$x_i$}}

\put(-35,40){\tiny{$x_j$}}

\put(-70,0){\tiny{$x_k$}}

\put(-35,0){\tiny{$x_l$}}


\multiput(-20,30)(7,0){6}{\line(1,0){5}}

\multiput(-20,10)(7,0){6}{\line(1,0){5}}

\put(20,30){\vector(1,0){0}}

\put(-20,10){\vector(-1,0){0}}


\linethickness{.5pt}

\put(47.5,10){\vector(1,1){0}}

\qbezier(40,0)(40,5)(47.5,10)

\qbezier(40,40)(40,35)(47.5,30)

\put(40,44){\vector(0,1){0}}

\put(52.5,10){\vector(-1,1){0}}

\qbezier(60,0)(60,5)(52.5,10)

\qbezier(60,40)(60,35)(52.5,30)

\put(60,44){\vector(0,1){0}}

\put(30,40){\tiny{$x_i$}}

\put(65,40){\tiny{$x_j$}}

\put(30,0){\tiny{$x_k$}}

\put(65,0){\tiny{$x_l$}}

\linethickness{5pt}

\put(50,10){\vector(0,1){20}}

\put(48,-15){\small{$\Gamma_1$}}

\end{picture}

\caption{$\Gamma_0$ and $\Gamma_1$}\label{maps}

\end{figure}
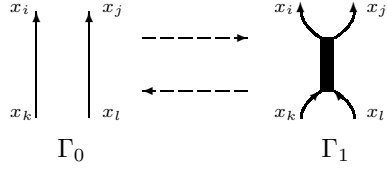

Let $\Gamma_0$ and $\Gamma_1$ be the two planar diagrams depicted
in Figure \ref{maps}. Then the matrix factorization
$C_n(\Gamma_0)$ is
\[
\left(%
\begin{array}{cc}
  \pi_{ik} & x_i-x_k \\
  \pi_{jl} & x_j-x_l \\
\end{array}%
\right)_{R'}.
\]
Explicitly, this is
\[
\left[%
\begin{array}{c}
  R' \\
  R'\{2-2n\} \\
\end{array}%
\right] \xrightarrow{P^{(n)}_0}
\left[%
\begin{array}{c}
  R'\{1-n\} \\
  R'\{1-n\} \\
\end{array}%
\right] \xrightarrow{P^{(n)}_1}
\left[%
\begin{array}{c}
  R' \\
  R'\{2-2n\} \\
\end{array}%
\right],
\]
where
\[
P^{(n)}_0=\left(%
\begin{array}{cc}
  \pi_{ik} & x_j-x_l \\
  \pi_{jl} & -x_i+x_k \\
\end{array}%
\right), \hspace{1cm}
P^{(n)}_1=\left(%
\begin{array}{cc}
  x_i-x_k & x_j-x_l \\
  \pi_{jl} & -\pi_{ik} \\
\end{array}%
\right).
\]
The matrix factorization $C_n(\Gamma_1)$ is
\[
\left(%
\begin{array}{cc}
  u_1 & x_i+x_j-x_k-x_l \\
  u_2 & x_ix_j-x_kx_l \\
\end{array}%
\right)_{R'}\{-1\}.
\]
Explicitly, this is
\[
\left[%
\begin{array}{c}
  R'\{-1\} \\
  R'\{3-2n\} \\
\end{array}%
\right] \xrightarrow{Q^{(n)}_0}
\left[%
\begin{array}{c}
  R'\{-n\} \\
  R'\{2-n\} \\
\end{array}%
\right] \xrightarrow{Q^{(n)}_1}
\left[%
\begin{array}{c}
  R'\{-1\} \\
  R'\{3-2n\} \\
\end{array}%
\right],
\]
where
\[
Q^{(n)}_0=\left(%
\begin{array}{cc}
  u_1 & x_ix_j-x_kx_l \\
  u_2 & -x_i-x_j+x_k+x_l \\
\end{array}%
\right), \hspace{.4cm}
Q^{(n)}_1=\left(%
\begin{array}{cc}
  x_i+x_j-x_k-x_l & x_ix_j-x_kx_l \\
  u_2 & -u_1 \\
\end{array}%
\right).
\]

Define $\chi^{(n)}_0:C_n(\Gamma_0)\rightarrow C_n(\Gamma_1)$ by
the matrices
\[
U^{(n)}_0=\left(%
\begin{array}{cc}
  x_k-x_j & 0 \\
  a_1 & 1 \\
\end{array}%
\right), \hspace{1cm} U^{(n)}_1=\left(%
\begin{array}{cc}
  x_k & -x_j \\
  -1 & 1 \\
\end{array}%
\right),
\]
and $\chi^{(n)}_1:C_n(\Gamma_1)\rightarrow C_n(\Gamma_0)$ by the
matrices
\[
V^{(n)}_0=\left(%
\begin{array}{cc}
  1 & 0 \\
  -a_1 & x_k-x_j \\
\end{array}%
\right), \hspace{1cm} V^{(n)}_1=\left(%
\begin{array}{cc}
  1 & x_j \\
  1 & x_k \\
\end{array}%
\right),
\]
where $a_1=-u_2+(u_1+x_iu_2-\pi_{jl})/(x_i-x_k)$.  These are
homomorphisms of matrix factorizations of degree $1$, i.e., these
commute with the differential maps and raise the quantum grading by $1$. So
$\chi^{(n)}_0$ and $\chi^{(n)}_1$ induce homomorphisms of
cohomologies. One can check that
\[
\chi^{(n)}_1\chi^{(n)}_0=(x_k-x_j)\id_{C_n(\Gamma_0)},
\hspace{1cm}
\chi^{(n)}_0\chi^{(n)}_1=(x_k-x_j)\id_{C_n(\Gamma_1)}.
\]

There is a similar construction for $C(\Gamma_0)$ and
$C(\Gamma_1)$. Explicitly, the matrix factorization $C(\Gamma_0)$ is
\[
\left[%
\begin{array}{c}
  R \\
  R\{-2,2\} \\
\end{array}%
\right] \xrightarrow{P_0}
\left[%
\begin{array}{c}
  R\{-1,1\} \\
  R\{-1,1\} \\
\end{array}%
\right] \xrightarrow{P_1}
\left[%
\begin{array}{c}
  R \\
  R\{-2,2\} \\
\end{array}%
\right],
\]
where
\[
P_0=\left(%
\begin{array}{cc}
  a & x_j-x_l \\
  a & -x_i+x_k \\
\end{array}%
\right), \hspace{1cm}
P_1=\left(%
\begin{array}{cc}
  x_i-x_k & x_j-x_l \\
  a & -a \\
\end{array}%
\right),
\]
and the matrix factorization $C(\Gamma_1)$ is
\[
\left[%
\begin{array}{c}
  R \\
  R\{-2,4\} \\
\end{array}%
\right] \xrightarrow{Q_0}
\left[%
\begin{array}{c}
  R\{-1,1\} \\
  R\{-1,3\} \\
\end{array}%
\right] \xrightarrow{Q_1}
\left[%
\begin{array}{c}
  R \\
  R\{-2,4\} \\
\end{array}%
\right],
\]
where
\[
Q_0=\left(%
\begin{array}{cc}
  a & x_ix_j-x_kx_l \\
  0 & -x_i-x_j+x_k+x_l \\
\end{array}%
\right), \hspace{.4cm}
Q_1=\left(%
\begin{array}{cc}
  x_i+x_j-x_k-x_l & x_ix_j-x_kx_l \\
  0 & -a \\
\end{array}%
\right).
\]

Define $\chi_0:C(\Gamma_0)\rightarrow C(\Gamma_1)$ by the matrices
\[
U_0=\left(%
\begin{array}{cc}
  x_k-x_j & 0 \\
  0 & 1 \\
\end{array}%
\right), \hspace{1cm} U_1=\left(%
\begin{array}{cc}
  x_k & -x_j \\
  -1 & 1 \\
\end{array}%
\right),
\]
and $\chi_1:C(\Gamma_1)\rightarrow C(\Gamma_0)$ by the matrices
\[
V_0=\left(%
\begin{array}{cc}
  1 & 0 \\
  0 & x_k-x_j \\
\end{array}%
\right), \hspace{1cm} V^{(n)}_1=\left(%
\begin{array}{cc}
  1 & x_j \\
  1 & x_k \\
\end{array}%
\right).
\]
These are homomorphisms of matrix factorizations, and induce
homomorphisms of cohomologies. One can check that $\chi_0$ and
$\chi_1$ are of bi-degrees $(0,2)$ and $(0,0)$, respectively, and
\[
\chi_1\chi_0=(x_k-x_j)\id_{C(\Gamma_0)}, \hspace{1cm}
\chi_0\chi_1=(x_k-x_j)\id_{C(\Gamma_1)}.
\]

\begin{figure}[h]

\setlength{\unitlength}{1pt}

\begin{picture}(420,120)(-210,-30)


\put(-80,10){\vector(0,1){40}}

\put(-100,10){\vector(0,1){40}}

\put(-95,0){$P_0$}


\put(-20,-30){\line(1,1){15}}

\put(5,-5){\vector(1,1){15}}

\put(20,-30){\vector(-1,1){40}}

\put(-5,-30){$P_-$}

\put(20,50){\line(-1,1){15}}

\put(-5,75){\vector(-1,1){15}}

\put(-20,50){\vector(1,1){40}}

\put(-5,50){$P_+$}

\multiput(30,70)(7,-3.5){6}{\line(2,-1){2.5}}

\put(70,50){\vector(2,-1){0}}

\put(50,65){\small{$1$}}

\multiput(-30,70)(-7,-3.5){6}{\line(-2,-1){2.5}}

\put(-70,50){\vector(-2,-1){0}}

\put(-55,65){\small{$0$}}

\multiput(30,-10)(7,3.5){6}{\line(2,1){2.5}}

\put(70,10){\vector(2,1){0}}

\put(50,-10){\small{$1$}}

\multiput(-30,-10)(-7,3.5){6}{\line(-2,1){2.5}}

\put(-70,10){\vector(-2,1){0}}

\put(-55,-10){\small{$0$}}


\linethickness{.5pt}

\put(87.5,20){\vector(1,1){0}}

\qbezier(80,10)(80,15)(87.5,20)

\qbezier(80,50)(80,45)(87.5,40)

\put(80,54){\vector(0,1){0}}

\put(92.5,20){\vector(-1,1){0}}

\qbezier(100,10)(100,15)(92.5,20)

\qbezier(100,50)(100,45)(92.5,40)

\put(100,54){\vector(0,1){0}}

\linethickness{5pt}

\put(90,20){\vector(0,1){20}}

\put(85,0){$P_1$}

\end{picture}

\caption{Resolutions of a crossing}\label{resolutions}

\end{figure}
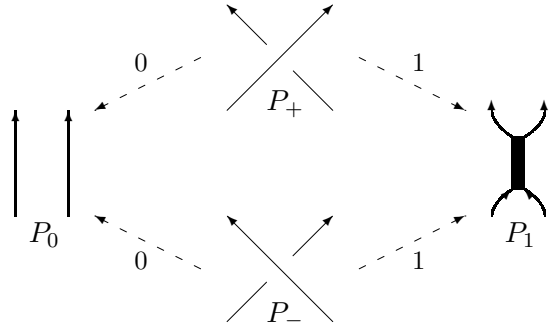

Let $D$ be an oriented link diagram. We put markings $\{x_1,\cdots,~x_m\}$
on $D$ so that none of the crossings is marked, and each arc
between two crossings has a marked point. For each crossing in
$D$, there are two ways to resolve it as shown in Figure
\ref{resolutions}. A resolution of $D$ is a planar diagram obtain
from $D$ by resolving all the crossings of $D$. Together with the
markings inherited from $D$, each resolution of $D$ is a planar
diagram satisfying the conditions in Subsection \ref{pd}.

For an arc $L^i_j$ from $x_j$ to $x_i$ that contains no crossings
and no other marked points, define $C_n(L^i_j)$ as above, and
consider it as the chain complex
\[
0 \rightarrow \underbrace{C_n(L^i_j)}_{0} \rightarrow 0,
\]
where $C_n(L^i_j)$ has cohomological degree $0$.

For a positive crossing $P_+$ in $D$, let $P_0$ and $P_1$ be the two
diagrams depicted in Figure \ref{resolutions}. Define $C_n(P_+)$
to be the chain complex
\[
0 \rightarrow \underbrace{C_n(P_1)\{n\}}_{-1} \xrightarrow{\chi^{(n)}_1}
\underbrace{C_n(P_0)\{n-1\}}_{0} \rightarrow 0,
\]
where $C_n(P_0)\{n-1\}$ has cohomological degree $0$, and
$C_n(P_1)\{n\}$ has cohomological degree $-1$.

For a negative crossing $P_-$ in $D$, define $C_n(P_+)$ to be the
chain complex
\[
0 \rightarrow \underbrace{C_n(P_0)\{1-n\}}_{0} \xrightarrow{\chi^{(n)}_0}
\underbrace{C_n(P_1)\{-n\}}_{1} \rightarrow 0,
\]
where $C_n(P_0)\{1-n\}$ has cohomological degree $0$, and
$C_n(P_1)\{-n\}$ has cohomological degree $1$.

Now define $C_n(D)$ to be the chain complex
\[
C_n(D) = (\bigotimes_{L^i_j}C_n(L^i_j)) ~\bigotimes~
(\bigotimes_{P}C_n(P)),
\]
where $L^i_j$ runs through all arcs in $D$ starting and ending in
marked points that contain no crossings and no other marked
points, and $P$ runs through all the crossings of $D$.

\begin{remark}
\begin{enumerate}
    \item Disregarding the chain complex structure on $C_n(D)$, we
    have
    \[
    C_n(D)=\bigoplus_{\Gamma}C_n(\Gamma)\{p_\Gamma\},
    \]
    where $\Gamma$ runs through all resolutions of $D$,
    $p_\Gamma\in\zed$. For a given $\Gamma$, all elements of
    $C_n(\Gamma)$ have the same cohomological degree.
    \item There are two differential maps on $C_n(D)$. One comes from
    the matrix factorizations $C_n(\Gamma)$, which we denote by
    $d_{mf}$. The other comes from the above construction, which we denote by
    $d_\chi$. These two differentials commute.
    \item $C_n(D)$ is $\zed\oplus\zed\oplus\zed_2$-graded, where
    the first $\zed$-grading is the cohomological grading, the
    second is inherited from $\zed$-grading of
    $C_n(\Gamma)\{p_\Gamma\}$, and is called the quantum grading,
    and the $\zed_2$-grading comes from the cyclic chain complex structure of
    $C_n(\Gamma)$. Note that $d_\chi$ preserves the quantum
    grading.
\end{enumerate}
\end{remark}

$d_\chi$ induces a differential map on the cohomology
$H(C_n(D),d_{mf})$, which we still denote by $d_\chi$.
Define $H_n(D)$ to be the cohomology
\[
H_n(D)=H(H(C_n(D),d_{mf}),d_\chi).
\]
$H_n(D)$ inherits the $\zed\oplus\zed$-gradings from $C_n(D)$. The
induced $\zed_2$-grading on $H_n(D)$ is trivial since all the
non-zero elements concentrate on the degree equal to the
$\zed_2$-number of components of the link represented by $D$.

M. Khovanov and L. Rozansky proved in \cite{KR1} that the
$\zed\oplus\zed$-graded isomorphism type of $H_n(D)$ is invariant
under Reidemeister moves, and, therefore, is a link invariant. They
also showed that, up to multiplication by a non-zero scalar, the
homomorphism of $H_n$ induced by a link cobodism is independent of the
movie representation. Specially this implies that the $H_n$ theory
gives bona fide homology groups, not just isomorphism classes.

The homology $H(D)$ is defined similarly. To $L^i_j$, $P_+$ and
$P_-$, we assign the chain complexes
\[
C(L^i_j)= [0 \rightarrow \underbrace{C(L^i_j)}_{0} \rightarrow 0],
\]
\[
C(P_+)= [0 \rightarrow \underbrace{C(P_1)\{0,-2\}}_{0}
\xrightarrow{\chi_1} \underbrace{C(P_0)\{0,-2\}}_{1} \rightarrow 0],
\]
\[
C(P_-)= [0 \rightarrow \underbrace{C(P_0)\{0,2\}}_{-1}
\xrightarrow{\chi_0} \underbrace{C(P_1)}_{0} \rightarrow 0],
\]
where the number below each term indicates its cohomological
degree. Then define
\[
C(D) = (\bigotimes_{L^i_j}C(L^i_j)) ~\bigotimes~ (\bigotimes_P
C(P)),
\]
where $L^i_j$ runs through all arcs in $D$ starting and ending in
marked points that contain no crossings and no other marked
points, and $P$ runs through all the crossings of $D$. Again, we
have
\[
C(D)=\bigoplus_{\Gamma}C(\Gamma)\{0, q_\Gamma\},
\]
where $\Gamma$ runs through all resolutions of $D$,
$q_\Gamma\in\zed$. Again, there are two differential maps:
$d_{mf}$ from the matrix factorization structure, and $d_\chi$
from the above chain complex construction. These commute with each
other. And $d_\chi$ preserves the two quantum gradings. We define
\[
H(D)=H(H(C(D),d_{mf}),d_\chi).
\]
$H(D)$ inherits the $\zed\oplus\zed\oplus\zed$-grading of $C(D)$,
i.e. the cohomological grading, the quantum $a$-grading and the
quantum $x$-grading. M. Khovanov and L. Rozansky proved in
\cite{KR1} that the $\zed\oplus\zed\oplus\zed$-graded isomorphism
type of $H(D)$ is invariant under Markov moves of braids up to
over all shift of grading. Thus, any two braid that represent the
same link have isomorphic $H$-invariant. Given a braid $B$ with
$b$ strands and writhe $w$, define
\[
\hat{H}(B) = H(B)\{\frac{w+b}{2},\frac{w+b}{2}\}[-\frac{w+b}{2}],
\]
where $[-\frac{w+b}{2}]$ means shifting the homological grading by
$-\frac{w+b}{2}$. From the construction in \cite{KR2}, it's easy to
see that the
$\frac{1}{2}\zed\oplus\frac{1}{2}\zed\oplus\frac{1}{2}\zed$-graded
isomorphism type of $\hat{H}(B)$ is honestly invariant under
Markov moves of braids.

\subsection{Graded Euler characteristics}

We normalize the HOMFLY polynomial $P$ by the skein relation
\[
\left\{
  \begin{array}{l}
    xP(L_-)-x^{-1}P(L_+)= yP(L_0), \\
    P(\text{unknot})=\frac{x-x^{-1}}{y},
  \end{array}
\right.
\]
where the diagrams of $L_+$, $L_-$ and $L_0$ are identical except at
one crossing, where they look like $P_+$, $P_-$ and $P_0$ in Figure
\ref{resolutions}, respectively.

For a link $L$, consider the graded Euler characteristic
\[
F_n(L)=\sum_{j,k}(-1)^jt^k\dim{H_n^{j,k}(L)}.
\]
From \cite{KR1}, we know that $F_n$ is characterized by the
following skein relation:
\[
\left\{
  \begin{array}{l}
    t^nF_n(L_-)-t^{-n}F_n(L_+)= (t-t^{-1})F_n(L_0), \\
    F_n(\text{unknot})=\frac{t^n-t^{-n}}{t-t^{-1}}.
  \end{array}
\right.
\]
So
\[
F_n(L)(t)=P(L)(t^n,t-t^{-1}).
\]

Next, consider the graded Euler characteristic
\[
\hat{F}(L)=\sum_{j,k,l}(-1)^jt^kq^l\dim{\hat{H}^{j,k,l}(L)}.
\]
From \cite{KR2}, we know that $\hat{F}$ is characterized by the
following skein relation:
\[
\left\{
  \begin{array}{l}
    -q^{-1}(-tq)^{\frac{1}{2}}\hat{F}(L_-)+q(-tq)^{-\frac{1}{2}}\hat{F}(L_+)= (q-q^{-1})\hat{F}(L_0), \\
    \hat{F}(\text{unknot})=-\frac{\sqrt{-1}t^{-\frac{1}{2}}q^{\frac{1}{2}}}{q-q^{-1}}.
  \end{array}
\right.
\]
where $L_+$, $L_-$ and $L_0$ are represented by braid diagrams that
identical except at one crossing, where they look like $P_+$, $P_-$
and $P_0$ in Figure \ref{resolutions}, respectively.

Change the variables by
\[
\left\{%
\begin{array}{ll}
    x = -q^{-1}(-tq)^{\frac{1}{2}} \\
    y = q-q^{-1} \\
\end{array}%
\right.
\]
Then, in these new variables, we have
\[
\left\{
  \begin{array}{l}
    x\hat{F}(L_-)-x^{-1}\hat{F}(L_+)= y\hat{F}(L_0), \\
    \hat{F}(\text{unknot})=-x^{-1}y^{-1}.
  \end{array}
\right.
\]

Compare these normalizations with those used in \cite{FW,Mo}. It's
clear that Theorem \ref{main} and Corollary \ref{braid-number} imply
the corresponding results in \cite{FW,Mo}. Specially, the two upper
bounds of the self-linking number in Theorem \ref{tb} are sharper
than the bound given by the HOMFLY polynomial, which is implicitly
proven in \cite{FW,Mo} (c.f. \cite{Fer,FT}).

\section{Resolved Braids}\label{inequality}

For positive integers $b$, $i$ with $1\leq i\leq b-1$, let
$\tau_i$ be the diagram depicted in Figure \ref{ti}. That is, from
left to right, $\tau_i$ consists of $i-1$ upward vertical regular
edges, then a vertical wide edge with two regular edges entering
through the bottom and two regular edges exiting through the top,
and then $b-i-1$ more upward vertical regular edges.

\begin{figure}[h]

\setlength{\unitlength}{1pt}

\begin{picture}(420,45)(-210,-10)

\linethickness{.5pt}

\put(-70,0){\vector(0,1){34}}

\put(-55,15){$\cdots$}

\put(-30,0){\vector(0,1){34}}

\put(-2.5,10){\vector(1,1){0}}

\qbezier(-10,0)(-10,5)(-2.5,10)

\qbezier(-10,30)(-10,25)(-2.5,20)

\put(-10,34){\vector(0,1){0}}

\put(2.5,10){\vector(-1,1){0}}

\qbezier(10,0)(10,5)(2.5,10)

\qbezier(10,30)(10,25)(2.5,20)

\put(10,34){\vector(0,1){0}}

\put(70,0){\vector(0,1){34}}

\put(45,15){$\cdots$}

\put(30,0){\vector(0,1){34}}

\put(-72,-10){\tiny{$1$}}

\put(-40,-10){\tiny{$i-1$}}

\put(-12,-10){\tiny{$i$}}

\put(2,-10){\tiny{$i+1$}}

\put(25,-10){\tiny{$i+2$}}

\put(67,-10){\tiny{$b$}}

\linethickness{5pt}

\put(0,10){\vector(0,1){10}}

\end{picture}

\caption{$\tau_i$}\label{ti}

\end{figure}

\begin{figure}[h]

\setlength{\unitlength}{1pt}

\begin{picture}(420,200)(-210,-60)


\linethickness{.5pt}

\put(-210,0){\vector(0,1){60}}

\put(-190,0){\vector(0,1){30}}

\put(-162.5,10){\vector(1,1){0}}

\qbezier(-170,0)(-170,5)(-162.5,10)

\qbezier(-170,30)(-170,25)(-162.5,20)

\put(-170,30){\vector(-1,4){0}}

\put(-157.5,10){\vector(-1,1){0}}

\qbezier(-150,0)(-150,5)(-157.5,10)

\qbezier(-150,30)(-150,25)(-157.5,20)

\qbezier(-190,30)(-190,35)(-182.5,40)

\qbezier(-190,60)(-190,55)(-182.5,50)

\put(-190,60){\vector(-1,4){0}}

\qbezier(-170,30)(-170,35)(-177.5,40)

\qbezier(-170,60)(-170,55)(-177.5,50)

\put(-150,30){\vector(0,1){60}}

\qbezier(-210,60)(-210,65)(-202.5,70)

\qbezier(-210,90)(-210,85)(-202.5,80)

\put(-210,90){\vector(-1,4){0}}

\qbezier(-190,60)(-190,65)(-197.5,70)

\qbezier(-190,90)(-190,85)(-197.5,80)

\put(-190,90){\vector(1,4){0}}

\put(-170,60){\vector(0,1){30}}

\put(-190,-60){$\tau_3\tau_2\tau_1$}

\linethickness{5pt}

\put(-160,10){\line(0,1){10}}

\put(-180,40){\line(0,1){10}}

\put(-200,70){\line(0,1){10}}


\linethickness{.5pt}

\put(0,0){\vector(0,1){60}}

\put(20,0){\vector(0,1){30}}

\qbezier(40,0)(40,5)(47.5,10)

\qbezier(40,30)(40,25)(47.5,20)

\put(40,30){\vector(-1,4){0}}

\qbezier(60,0)(60,5)(52.5,10)

\qbezier(60,30)(60,25)(52.5,20)

\qbezier(20,30)(20,35)(27.5,40)

\qbezier(20,60)(20,55)(27.5,50)

\put(20,60){\vector(-1,4){0}}

\qbezier(40,30)(40,35)(32.5,40)

\qbezier(40,60)(40,55)(32.5,50)

\put(60,30){\vector(0,1){60}}

\qbezier(0,60)(0,65)(7.5,70)

\qbezier(0,90)(0,85)(7.5,80)

\put(0,90){\vector(-1,4){0}}

\qbezier(20,60)(20,65)(12.5,70)

\qbezier(20,90)(20,85)(12.5,80)

\put(20,90){\vector(1,4){0}}

\put(40,60){\vector(0,1){30}}

\put(70,-60){$\underline{\tau_3\tau_2\tau_1}$}

\put(100,90){\vector(0,-10){90}}

\put(120,90){\vector(0,-10){90}}

\put(140,90){\vector(0,-10){90}}

\put(160,90){\vector(0,-10){90}}

\qbezier(60,90)(60,100)(80,100)

\qbezier(40,90)(40,110)(80,110)

\qbezier(20,90)(20,120)(80,120)

\qbezier(0,90)(0,130)(80,130)

\qbezier(100,90)(100,100)(80,100)

\qbezier(120,90)(120,110)(80,110)

\qbezier(140,90)(140,120)(80,120)

\qbezier(160,90)(160,130)(80,130)

\qbezier(60,0)(60,-10)(80,-10)

\qbezier(40,0)(40,-20)(80,-20)

\qbezier(20,0)(20,-30)(80,-30)

\qbezier(0,0)(0,-40)(80,-40)

\qbezier(100,0)(100,-10)(80,-10)

\qbezier(120,0)(120,-20)(80,-20)

\qbezier(140,0)(140,-30)(80,-30)

\qbezier(160,0)(160,-40)(80,-40)

\linethickness{5pt}

\put(50,10){\line(0,1){10}}

\put(30,40){\line(0,1){10}}

\put(10,70){\line(0,1){10}}

\end{picture}

\caption{$\tau_3\tau_2\tau_1$ and $\underline{\tau_3\tau_2\tau_1}$
for $b=4$.}\label{eg}

\end{figure}
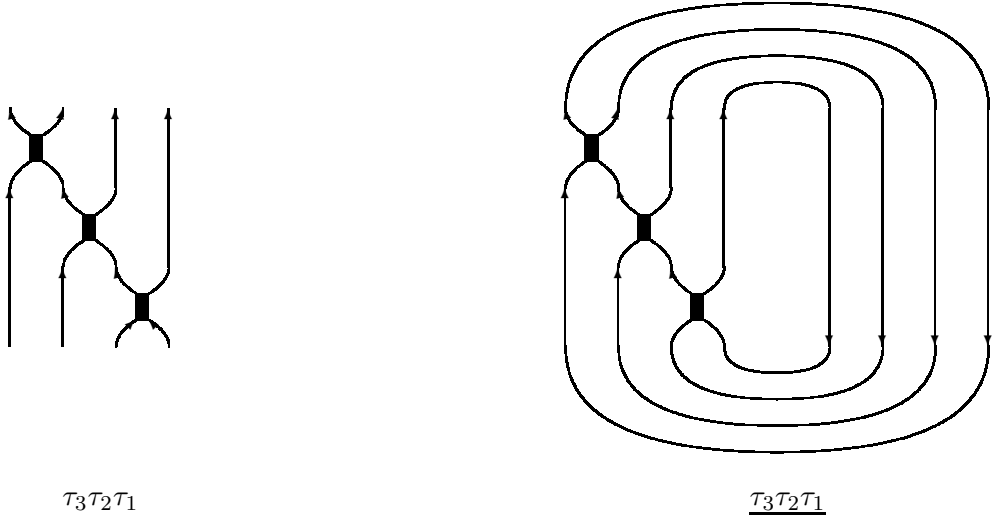

We use the word $\tau_{i_1}\cdots\tau_{i_m}$ to represent the planar
graph formed by stacking the graphs $\tau_{i_1}$, $\cdots$,
$\tau_{i_m}$ together vertically from bottom to top with the top end
points of $\tau_{i_l}$ identified with the corresponding bottom end
points of $\tau_{i_{l+1}}$. The symbol
$\underline{\tau_{i_1}\cdots\tau_{i_m}}$ represents the closed graph
obtained from $\tau_{i_1}\cdots\tau_{i_m}$ by attaching a disjoint
regular edge from each end point on the top to the corresponding end
point at the bottom. We use the convention that the empty word
$\phi$ represents $b$ vertical upward regular edges, and, therefore,
$\underline{\phi}$ represents $b$ concentric circles. We call
$\tau_{i_1}\cdots\tau_{i_m}$ a resolved braid, and
$\underline{\tau_{i_1}\cdots\tau_{i_m}}$ a resolved closed braid.
See Figure \ref{eg} for examples of such graphs. There are two obvious types of isotopies:

(I) If $|i-j|>1$, then $\tau_i\tau_j$ is isotopic to
$\tau_j\tau_i$;

(II) If $\mu$ and $\nu$ are two words in $\tau_{1}$, $\cdots$,
$\tau_{b-1}$, then $\underline{\mu\nu}$ is isotopic to
$\underline{\nu\mu}$.

We prove inequality \eqref{main-inequality} in Theorem \ref{main}
first.

For a word $\mu$ in $\tau_{1}$, $\cdots$, $\tau_{b-1}$, consider the
$\mathbb{Z}\oplus\mathbb{Z}$-graded cohomology
$H(\underline{\mu})=H(C(\underline{\mu}),d_{mf})$. Denote by
$H^{k,l}(\underline{\mu})$ the subspace of $H(\underline{\mu})$
consisting of homogeneous elements of $\mathbb{Z}\oplus\mathbb{Z}$
grading $(k,l)$, where $k$ is the quantum $a$-grading, and $l$ is
the quantum $x$-grading. Define
\[
g_{max}(\underline{\mu})=\max\{k~|~H^{k,l}(\underline{\mu})\neq0
\text{ for some }l\},
\]
\[
g_{min}(\underline{\mu})=\min\{k~|~H^{k,l}(\underline{\mu})\neq0
\text{ for some }l\}.
\]

In order to calculate $H(\underline{\mu})$, we need put marks on
regular edges of $\underline{\mu}$. But, as long as there is at
least one mark on each regular edge, the result is independent of
the choice of the marks.

\begin{lemma}\label{circles}
For the empty word $\phi$ (with $b$ strands), we have
$g_{max}(\underline{\phi})=-1$ and $g_{min}(\underline{\phi})=-b$.
\end{lemma}
\begin{proof}
$\underline{\phi}$ represents $b$ concentric circles. We order the
circles, and put the marking $x_i$ on the $i$-th circle. The cyclic
Koszul complex $C(\underline{\phi})$ is given over
$R=\mathbb{Q}[a,x_1,\cdots,x_b]$ by
\[
\left(%
\begin{array}{cc}
  a & 0 \\
  \vdots & \vdots \\
  a & 0 \\
\end{array}%
\right),
\]
where there are $b$ rows. After $b-1$ elementary transformations of
rows(c.f. subsection 2.1 of \cite{KR2}), we see that it is isomorphic
to
\[
\left(%
\begin{array}{cc}
  a & 0 \\
  0 & 0 \\
  \vdots & \vdots \\
  0 & 0 \\
\end{array}%
\right),
\]
which is isomorphic to the chain complex
\[
C \otimes_{\mathbb{Q}} D^{\otimes(b-1)} \otimes_{\mathbb{Q}} \mathbb{Q}[x_1,\cdots,x_b],
\]
where $C$ is the chain complex 
\[
0 \rightarrow \mathbb{Q}[a] \stackrel{a}{\rightarrow} \mathbb{Q}[a] \{-1,1\} \rightarrow 0,
\]
and $D$ is the chain complex
\[
0 \rightarrow \mathbb{Q} \stackrel{0}{\rightarrow} \mathbb{Q} \{-1,1\} \rightarrow 0.
\]
Then, by the K$\ddot{\text{u}}$nneth Theorem and the Universal Coefficient Theorem, it's easy to check that 
\begin{eqnarray*}
H(\underline{\phi}) & = & H(C \otimes_{\mathbb{Q}} D^{\otimes(b-1)} \otimes_{\mathbb{Q}} \mathbb{Q}[x_1,\cdots,x_b]) \\
& = & H(C) \otimes_{\mathbb{Q}} H(D^{\otimes(b-1)}) \otimes_{\mathbb{Q}} \mathbb{Q}[x_1,\cdots,x_b] \\
& = & \mathbb{Q}\{-1,1\} \otimes_{\mathbb{Q}} D^{\otimes(b-1)} \otimes_{\mathbb{Q}} \mathbb{Q}[x_1,\cdots,x_b].
\end{eqnarray*}
From the definition of $D$, one can see that the maximal and minimal quantum $a$-gradings of the right hand side are $-1$ and $-b$. This proves the lemma.
\end{proof}

\begin{lemma}\label{move1}
If $\mu=\tau_{i_1}\cdots\tau_{i_m}$ is a word in $\tau_{1}$,
$\cdots$, $\tau_{b-1}$ with $i_p<i$ for $1\leq p \leq m$, and
$\mu'=\mu\tau_i$, then
$g_{max}(\underline{\mu})=g_{max}(\underline{\mu'})$, and
$g_{min}(\underline{\mu})=g_{min}(\underline{\mu'})$.
\end{lemma}
\begin{proof}
The only difference between $\underline{\mu}$ and $\underline{\mu'}$
occurs in the part depicted in Figure \ref{reduction1}.

\begin{figure}[h]

\setlength{\unitlength}{1pt}

\begin{picture}(420,60)(-210,-20)

\linethickness{.5pt}


\put(-90,-10){\vector(0,1){50}}

\put(-70,0){\line(0,1){30}}

\qbezier(-70,30)(-70,40)(-60,40)

\qbezier(-60,40)(-50, 40)(-50,30)

\put(-50,30){\vector(0,-1){30}}

\qbezier(-70,0)(-70,-10)(-60,-10)

\qbezier(-60,-10)(-50, -10)(-50,0)

\put(-75,-20){$\underline{\mu}$}

\put(-100,-10){\tiny{$x_4$}}

\put(-100,40){\tiny{$x_1$}}

\put(-45,15){\tiny{$x_2$}}

\put(-51.5,15){-}


\put(57.5,10){\vector(1,1){0}}

\qbezier(50,-10)(50,5)(57.5,10)

\qbezier(50,40)(50,25)(57,20)

\put(50,44){\vector(0,1){0}}

\qbezier(70,0)(70,5)(62.5,10)

\qbezier(70,30)(70,25)(62.5,20)

\qbezier(70,30)(70,40)(80,40)

\qbezier(80,40)(90, 40)(90,30)

\put(90,30){\vector(0,-1){30}}

\qbezier(70,0)(70,-10)(80,-10)

\qbezier(80,-10)(90, -10)(90,0)

\put(65,-20){$\underline{\mu'}$}

\put(40,-10){\tiny{$x_4$}}

\put(40,40){\tiny{$x_1$}}

\put(95,15){\tiny{$x_2$}}

\put(88.5,15){-}

\linethickness{5pt}

\put(60,10){\line(0,1){10}}

\end{picture}

\caption{$\underline{\mu}$ and
$\underline{\mu'}$}\label{reduction1}

\end{figure}
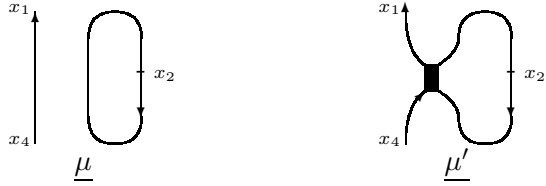

The local matrix factorization of the part of $\underline{\mu}$ in Figure \ref{reduction1} is
\[
\left(%
\begin{array}{cc}
  a & x_1-x_4 \\
  a & 0 \\
\end{array}%
\right).
\]
After an elementary transformation, it becomes
\[
\left(%
\begin{array}{cc}
  a & x_1-x_4 \\
  0 & 0 \\
\end{array}%
\right).
\]
Write $R=\mathbb{Q}[a,x_1,x_2,x_4]$. Explicitly, this matrix factorization is
\[
R \oplus R\{-2,2\} \rightarrow R\{-1,1\} \oplus R\{-1,1\} \rightarrow R \oplus R\{-2,2\},
\]
where
\[
d^0 =
\left(%
\begin{array}{cc}
  a & 0 \\
  0 & -x_1+x_4 \\
\end{array}%
\right)
~\text{ and }~
d^1 =
\left(%
\begin{array}{cc}
  x_1-x_4 & 0 \\
  0 & -a \\
\end{array}%
\right).
\]

The local matrix factorization of the part of $\underline{\mu'}$ in
Figure \ref{reduction1} is
\[
\left(%
\begin{array}{cc}
  a & x_1-x_4 \\
  0 & 0 \\
\end{array}%
\right).
\]

Explicitly, this matrix factorization is
\[
R \oplus R\{-2,4\} \rightarrow R\{-1,1\} \oplus R\{-1,3\}
\rightarrow R \oplus R\{-2,4\},
\]
where
\[
d^0 =
\left(%
\begin{array}{cc}
  a & 0 \\
  0 & -x_1+x_4 \\
\end{array}%
\right)
~\text{ and }~
d^1 =
\left(%
\begin{array}{cc}
  x_1-x_4 & 0 \\
  0 & -a \\
\end{array}%
\right).
\]

Let $f$ be the map between the two matrix factorizations given by
\[
\left(%
\begin{array}{cc}
  1 & 0 \\
  0 & 1 \\
\end{array}%
\right):
R \oplus R\{-2,2\} \rightarrow R \oplus R\{-2,4\},
\]
and
\[
\left(%
\begin{array}{cc}
  1 & 0 \\
  0 & 1 \\
\end{array}%
\right):
R\{-1,1\} \oplus R\{-1,1\} \rightarrow R\{-1,1\} \oplus R\{-1,3\}.
\]
It's clear that $f$ is an isomorphism of matrix factorizations that
preserves the quantum $a$-grading (but not the $x$-grading). After
tensoring $f$ with the identity maps of other local matrix
factorizations of $\underline{\mu}$ (and $\underline{\mu'}$), we get
an isomorphism between these matrix factorizations that preserves
the $a$-grading. This proves the lemma.
\end{proof}

\begin{lemma}\label{move2}
If $\nu$ is a word in $\tau_{1}$, $\cdots$, $\tau_{b-1}$,
$\mu=\nu\tau_i$ and $\mu'=\nu\tau_i\tau_i$, then
$g_{max}(\underline{\mu})=g_{max}(\underline{\mu'})$, and
$g_{min}(\underline{\mu})=g_{min}(\underline{\mu'})$.
\end{lemma}
\begin{proof}
(Follow Proposition 30 of \cite{KR1}.) The only difference between
$\underline{\mu}$ and $\underline{\mu'}$ occurs in the part depicted
in Figure \ref{reduction2}.

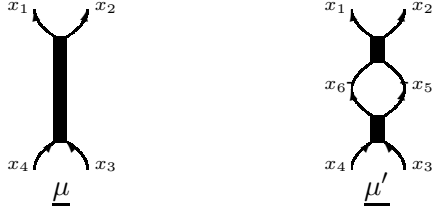
\begin{figure}[h]

\setlength{\unitlength}{1pt}

\begin{picture}(420,70)(-210,-10)


\linethickness{.5pt}

\put(62.5,10){\vector(-1,1){0}}

\qbezier(70,0)(70,5)(62.5,10)

\qbezier(70,30)(70,25)(62.5,20)

\put(70,30){\vector(1,4){0}}

\put(50,30){\vector(-1,4){0}}

\put(57.5,10){\vector(1,1){0}}

\qbezier(50,0)(50,5)(57.5,10)

\qbezier(50,30)(50,25)(57.5,20)

\qbezier(70,30)(70,35)(62.5,40)

\qbezier(70,60)(70,55)(62.5,50)

\put(70,60){\vector(1,4){0}}

\put(50,60){\vector(-1,4){0}}

\qbezier(50,30)(50,35)(57.5,40)

\qbezier(50,60)(50,55)(57.5,50)

\put(55,-10){$\underline{\mu'}$}

\put(40,60){\tiny{$x_1$}}

\put(40,0){\tiny{$x_4$}}

\put(73,60){\tiny{$x_2$}}

\put(73,0){\tiny{$x_3$}}

\put(40,30){\tiny{$x_6$}}

\put(73,30){\tiny{$x_5$}}

\put(48.5,30){-}

\put(68.5,30){-}

\linethickness{5pt}

\put(60,10){\line(0,1){10}}

\put(60,40){\line(0,1){10}}


\linethickness{.5pt}

\put(-62.5,10){\vector(1,1){0}}

\qbezier(-70,0)(-70,5)(-62.5,10)

\put(-57.5,10){\vector(-1,1){0}}

\qbezier(-50,0)(-50,5)(-57.5,10)

\qbezier(-70,60)(-70,55)(-62.5,50)

\put(-70,60){\vector(-1,4){0}}

\put(-50,60){\vector(1,4){0}}

\qbezier(-50,60)(-50,55)(-57.5,50)

\put(-63,-10){$\underline{\mu}$}

\put(-80,60){\tiny{$x_1$}}

\put(-80,0){\tiny{$x_4$}}

\put(-47,60){\tiny{$x_2$}}

\put(-47,0){\tiny{$x_3$}}

\linethickness{5pt}

\put(-60,10){\line(0,1){40}}

\end{picture}

\caption{$\underline{\mu}$ and
$\underline{\mu'}$}\label{reduction2}

\end{figure}

Write $R=\mathbb{Q}[a,x_1,\cdots,x_6]$. Let $s_1=x_5+x_6$,
$s_2=x_5x_6$, and $\check{R}=\mathbb{Q}[a,x_1,x_2,x_3,x_4,s_1,s_2]$.
Then $R=\check{R}\oplus \check{R}\{0,2\}$.

The local matrix factorization $M'$ of the part of
$\underline{\mu'}$ in Figure \ref{reduction1} is
\[
\left(%
\begin{array}{cc}
  a & x_1+x_2-s_1 \\
  0 & x_1x_2-s_2  \\
  a & s_1-x_3-x_4 \\
  0 & s_2-x_3x_4  \\
\end{array}%
\right)_R.
\]

Let $\check{M}'$ be the matrix factorization over $\check{R}$ given
by the same sequence. Then $M'=\check{M}'\oplus \check{M}'\{0,2\}$.
After elementary transformations, $\check{M}'$ is isomorphic to
\[
\left(%
\begin{array}{cc}
  a & x_1+x_2-x_3-x_4 \\
  0 & x_1x_2-x_3x_4  \\
  0 & s_1-x_3-x_4 \\
  0 & s_2-x_3x_4  \\
\end{array}%
\right)_{\check{R}}.
\]
We then exclude the variables $s_1$, $s_2$ and remove the last two
rows from the above sequence. This does not change chain homotopy
type of the matrix factorization (c.f. Proposition 3 of \cite{KR2}).
The result of this operation is a matrix factorization over
$\mathbb{Q}[a,x_1,x_2,x_3,x_4]$ given by the sequence
\[
\left(%
\begin{array}{cc}
  a & x_1+x_2-x_3-x_4 \\
  0 & x_1x_2-x_3x_4  \\
\end{array}%
\right),
\]
which is exactly the local matrix factorization of the part of
$\underline{\mu}$ depicted in Figure \ref{reduction2}. This implies
that $H(\underline{\mu'})=H(\underline{\mu})\oplus
H(\underline{\mu})\{0,2\}$. And the lemma follows.
\end{proof}

\begin{lemma}\label{move3}
If $\nu$ is a word in $\tau_{1}$, $\cdots$, $\tau_{b-1}$,
$\mu_1=\nu\tau_{i-1}\tau_i\tau_{i-1}$, $\mu_2=\nu\tau_i$, and
$\mu_3=\nu\tau_i\tau_{i-1}\tau_i$ then
$g_{max}(\underline{\mu_3})\leq
\max\{g_{max}(\underline{\mu_1}),g_{max}(\underline{\mu_2})\}$,
$g_{min}(\underline{\mu_3})\geq
\min\{g_{min}(\underline{\mu_1}),g_{min}(\underline{\mu_2})\}$.
\end{lemma}
\begin{proof}
Let $\mu_4=\nu\tau_{i-1}$. The only difference between
$\underline{\mu_1}$, $\underline{\mu_2}$, $\underline{\mu_3}$ and
$\underline{\mu_4}$ occurs in the part depicted in Figure
\ref{reduction3}.

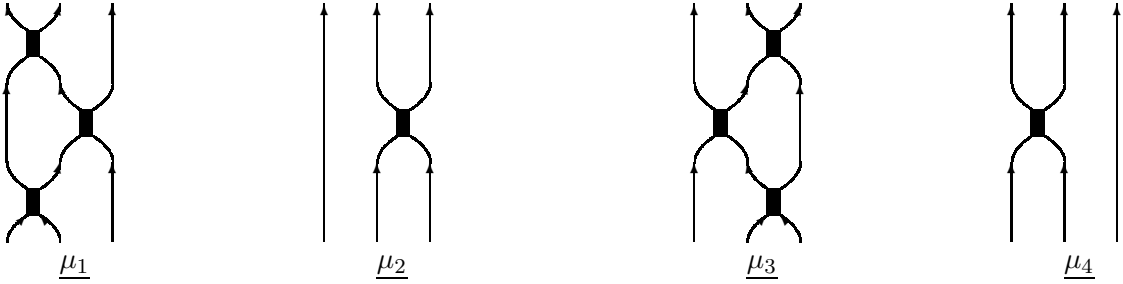
\begin{figure}[h]

\setlength{\unitlength}{1pt}

\begin{picture}(420,100)(-210,-10)


\linethickness{.5pt}

\put(-210,30){\vector(0,1){30}}

\put(-170,0){\vector(0,1){30}}

\put(-202.5,10){\vector(1,1){0}}

\qbezier(-210,0)(-210,5)(-202.5,10)

\qbezier(-210,30)(-210,25)(-202.5,20)

\put(-197.5,10){\vector(-1,1){0}}

\qbezier(-190,0)(-190,5)(-197.5,10)

\qbezier(-190,30)(-190,25)(-197.5,20)

\qbezier(-190,30)(-190,35)(-182.5,40)

\qbezier(-190,60)(-190,55)(-182.5,50)

\put(-190,60){\vector(-1,4){0}}

\qbezier(-170,30)(-170,35)(-177.5,40)

\qbezier(-170,60)(-170,55)(-177.5,50)

\qbezier(-210,60)(-210,65)(-202.5,70)

\qbezier(-210,90)(-210,85)(-202.5,80)

\put(-210,90){\vector(-1,4){0}}

\qbezier(-190,60)(-190,65)(-197.5,70)

\qbezier(-190,90)(-190,85)(-197.5,80)

\put(-190,90){\vector(1,4){0}}

\put(-190,30){\vector(1,4){0}}

\put(-170,60){\vector(0,1){30}}

\put(-190,-10){$\underline{\mu_1}$}

\linethickness{5pt}

\put(-200,10){\line(0,1){10}}

\put(-180,40){\line(0,1){10}}

\put(-200,70){\line(0,1){10}}


\linethickness{.5pt}

\put(-90,0){\vector(0,1){90}}

\put(-50,0){\vector(0,1){30}}

\put(-70,0){\vector(0,1){30}}

\qbezier(-50,30)(-50,35)(-57.5,40)

\qbezier(-50,60)(-50,55)(-57.5,50)

\qbezier(-70,30)(-70,35)(-62.5,40)

\qbezier(-70,60)(-70,55)(-62.5,50)

\put(-50,60){\vector(0,1){30}}

\put(-70,60){\vector(0,1){30}}

\put(-70,-10){$\underline{\mu_2}$}

\linethickness{5pt}

\put(-60,40){\line(0,1){10}}


\linethickness{.5pt}

\put(90,30){\vector(0,1){30}}

\put(50,0){\vector(0,1){30}}

\put(82.5,10){\vector(-1,1){0}}

\qbezier(90,0)(90,5)(82.5,10)

\qbezier(90,30)(90,25)(82.5,20)

\put(77.5,10){\vector(1,1){0}}

\qbezier(70,0)(70,5)(77.5,10)

\qbezier(70,30)(70,25)(77.5,20)

\qbezier(70,30)(70,35)(62.5,40)

\qbezier(70,60)(70,55)(62.5,50)

\put(70,60){\vector(1,4){0}}

\qbezier(50,30)(50,35)(57.5,40)

\qbezier(50,60)(50,55)(57.5,50)

\qbezier(90,60)(90,65)(82.5,70)

\qbezier(90,90)(90,85)(82.5,80)

\put(90,90){\vector(1,4){0}}

\qbezier(70,60)(70,65)(77.5,70)

\qbezier(70,90)(70,85)(77.5,80)

\put(70,90){\vector(-1,4){0}}

\put(70,30){\vector(-1,4){0}}

\put(50,60){\vector(0,1){30}}

\put(70,-10){$\underline{\mu_3}$}

\linethickness{5pt}

\put(80,10){\line(0,1){10}}

\put(60,40){\line(0,1){10}}

\put(80,70){\line(0,1){10}}


\linethickness{.5pt}

\put(210,0){\vector(0,1){90}}

\put(170,0){\vector(0,1){30}}

\put(190,0){\vector(0,1){30}}

\qbezier(170,30)(170,35)(177.5,40)

\qbezier(170,60)(170,55)(177.5,50)

\qbezier(190,30)(190,35)(182.5,40)

\qbezier(190,60)(190,55)(182.5,50)

\put(170,60){\vector(0,1){30}}

\put(190,60){\vector(0,1){30}}

\put(190,-10){$\underline{\mu_4}$}

\linethickness{5pt}

\put(180,40){\line(0,1){10}}

\end{picture}

\caption{$\underline{\mu_1}$, $\underline{\mu_2}$,
$\underline{\mu_3}$ and $\underline{\mu_4}$}\label{reduction3}

\end{figure}

By Proposition 7 of \cite{KR2},
\[
H(\underline{\mu_1}) \oplus H(\underline{\mu_2})\{0,2\} =
H(\underline{\mu_3}) \oplus H(\underline{\mu_4})\{0,2\}.
\]
This implies the lemma.
\end{proof}

\begin{lemma}\label{induction}
Let $\mu=\tau_{i_1}\cdots\tau_{i_m}$ be a word in $\tau_{1}$,
$\cdots$, $\tau_{b-1}$, such that $m\geq2$, $i_1=i_m=i$, and $i_p<i$
for $1<p<m$. Then, after possibly finitely many steps of isotopies
of type (I), there is a sub-word of $\mu$ of the form $\tau_j\tau_j$
or $\tau_j\tau_{j-1}\tau_j$ for some $j\leq i$.
\end{lemma}
\begin{proof}
We induct on the length $m$ of the word $\mu$. When $m=2$, the
lemma is trivially true. Assume that the lemma is true for words
with the given properties and length $=2,\cdots, m-1$. Now
consider $\mu=\tau_{i_1}\cdots\tau_{i_m}$. If there is none or
only one $\tau_{i-1}$ among $\tau_{i_2},\cdots,\tau_{i_{m-1}}$, then the
lemma is true for $\mu$ with $j=i$. If there are more than one
$\tau_{i-1}$'s among $\tau_{i_2},\cdots,\tau_{i_{m-1}}$, then there is a
sub-word $\nu=\tau_{j_1}\cdots\tau_{j_l}$ of $\mu$ with $2\leq
l<m$, $j_1=j_l=i-1$, and $j_q<i-1$ for $1<q<l$. Thus, by
induction, the lemma is true for $\mu$.
\end{proof}

\begin{proposition}\label{resolved-braids}
If $\mu$ is a word in $\tau_{1}$, $\cdots$, $\tau_{b-1}$, then
$-b\leq g_{min}(\underline{\mu})\leq g_{max}(\underline{\mu}) \leq
-1$.
\end{proposition}
\begin{proof}
For a word $\mu=\tau_{i_1}\cdots\tau_{i_m}$, define the weight
$W(\mu)$ of $\mu$ to be
\[
W(\mu)=\sum_{p=1}^m i_p.
\]

We prove the proposition by induction on the weight of a word. If
$W(\mu)=0$, then $\mu$ is the empty word $\phi$, and Lemma
\ref{circles} implies the proposition. Now assume that the
proposition is true for all words with weight $<W$, and
$\mu=\tau_{i_1}\cdots\tau_{i_m}$ is a word with weight $W$. Let
$i=\max\{i_1,\cdots,i_m\}$. If $i$ appears only once among
$i_1,\cdots,i_m$, then, after a possible type (II) isotopy, Lemma
\ref{move1} implies that the proposition is true for $\mu$. If $i$
appears at least twice among $i_1,\cdots,i_m$, then, after
possible type (I) and type (II) isotopies, Lemmas \ref{move2},
\ref{move3} and \ref{induction} imply that the proposition is true
for $\mu$.
\end{proof}

\begin{proof}[Proof of inequality \eqref{main-inequality} in Theorem
\ref{main}] For a braid $B$ with $b$ strands, consider the
Khovanov-Rozansky cohomology $H(B)$, and let
\[
g_{max}(B)=\max\{k~|~H^{j,k,l}(L)\neq0 \text{ for some }j,~l\},
\]
\[
g_{min}(B)=\min\{k~|~H^{j,k,l}(L)\neq0 \text{ for some }j,~l\}.
\]
Note that
\[
H(B)=H(H(C(B),d_{mf}),d_\chi),
\]
and
\[
H(C(B),d_{mf})=\bigoplus_{\Gamma}H(\Gamma)\{0, q_\Gamma\},
\]
where $\Gamma$ runs through all resolutions of $B$, and $q_\Gamma\in
\zed$. Thus, by Proposition \ref{resolved-braids}, we have $-b \leq
g_{min}(B) \leq g_{max}(B) \leq -1$. But
$\hat{H}(B)=H(B)\{\frac{w+b}{2},\frac{w+b}{2}\}[-\frac{w+b}{2}]$. So
\[
\hat{g}_{max}(B)=g_{max}(B)+\frac{w+b}{2},
\]
\[
\hat{g}_{min}(B)=g_{min}(B)+\frac{w+b}{2}.
\]
And inequality \eqref{main-inequality} follows.
\end{proof}

The proof of inequality \eqref{main-inequality-n} in Theorem
\ref{main} is similar.

For a word $\mu$ in $\tau_{1}$, $\cdots$, $\tau_{b-1}$, consider the
$\mathbb{Z}$-graded cohomology
$H_n(\underline{\mu})=H_n(C(\underline{\mu}),d_{mf})$. Denote by
$H_n^{k}(\underline{\mu})$ the subspace of $H_n(\underline{\mu})$
consisting of homogeneous elements of quantum grading $k$. Define
\[
g^{(n)}_{max}(\underline{\mu})=\max\{k~|~H_n^{k}(\underline{\mu})\neq0\},
\]
\[
g^{(n)}_{min}(\underline{\mu})=\min\{k~|~H_n^{k}(\underline{\mu})\neq0\}.
\]

We have following conclusions similar to Lemmas \ref{circles},
\ref{move1}, \ref{move2} and \ref{move3}.

\begin{lemma}\label{hn-moves}
\begin{enumerate}
\item Let $\phi$ be the empty word (with $b$ strands). Then
\[
g^{(n)}_{max}(\underline{\phi})=(n-1)b, ~\text{ and }~
g^{(n)}_{min}(\underline{\phi})=-(n-1)b.
\]
\item If $\mu=\tau_{i_1}\cdots\tau_{i_m}$ is a word in $\tau_{1}$,
$\cdots$, $\tau_{b-1}$ with $i_p<i$ for $1\leq p \leq m$, and
$\mu'=\mu\tau_i$, then
\[
g^{(n)}_{max}(\underline{\mu'})=g^{(n)}_{max}(\underline{\mu})-1,
~\text{ and }~
g^{(n)}_{min}(\underline{\mu'})=g^{(n)}_{min}(\underline{\mu})+1.
\]
\item If $\nu$ is a word in $\tau_{1}$, $\cdots$, $\tau_{b-1}$,
$\mu=\nu\tau_i$ and $\mu'=\nu\tau_i\tau_i$, then
\[
g^{(n)}_{max}(\underline{\mu'})=g^{(n)}_{max}(\underline{\mu})+1,
~\text{ and }~
g^{(n)}_{min}(\underline{\mu'})=g^{(n)}_{min}(\underline{\mu})-1.
\]
\item If $\nu$ is a word in $\tau_{1}$, $\cdots$, $\tau_{b-1}$,
$\mu_1=\nu\tau_{i-1}\tau_i\tau_{i-1}$, $\mu_2=\nu\tau_i$, and
$\mu_3=\nu\tau_i\tau_{i-1}\tau_i$ then
\[
g^{(n)}_{max}(\underline{\mu_3})\leq
\max\{g^{(n)}_{max}(\underline{\mu_1}),g^{(n)}_{max}(\underline{\mu_2})\},
~\text{ and }~ g^{(n)}_{min}(\underline{\mu_3})\geq
\min\{g^{(n)}_{min}(\underline{\mu_1}),g^{(n)}_{min}(\underline{\mu_2})\}.
\]
\end{enumerate}
\end{lemma}

\begin{proof}
(1) We mark the $i$-th strand by a single variable $x_i$ of degree $2$, for $i=1,\cdots, b$. A straightforward computation shows that
\[
H_n(\underline{\phi}) = A_1 \otimes_{\mathbb{Q}} A_2
\otimes_{\mathbb{Q}} \cdots \otimes_{\mathbb{Q}} A_b,
\]
where $A_i=[\mathbb{Q}[x_i]/(x_i^n)]\{-n+1\}$. The first part of the lemma follows from this.

(2) Note that the $(i+1)$-th strand in $\underline{\mu}$ is a circle. We mark it by a single variable $x$ of degree $2$. Denote by $\hat{\underline{\mu}}$ the closed resolved braid obtained from $\underline{\mu}$ by removing the $(i+1)$-th strand. Then
\[
H_n(\underline{\mu}) = H_n(\hat{\underline{\mu}})\otimes_{\mathbb{Q}}A = \bigoplus_{i=0}^{n-1}H_n(\hat{\underline{\mu}})\{-n+1+2i\},
\]
where $A=[\mathbb{Q}[x]/(x^n)]\{-n+1\}$. But, from Proposition 29 of \cite{KR1}, we know that
\[
H_n(\underline{\mu'}) = \bigoplus_{i=0}^{n-2}H_n(\hat{\underline{\mu}})\{-n+2+2i\}.
\]
Compare these two decompositions. And the second part of the lemma follows.

(3) By Proposition 30 of \cite{KR1}, we have
\[
H_n(\underline{\mu'})=H_n(\underline{\mu})\{1\}\oplus
H_n(\underline{\mu})\{-1\},
\]
which implies the third part of the lemma.

(4) By Proposition 32 of \cite{KR1}, we have
\[
H_n(\underline{\mu_1}) \oplus H_n(\underline{\mu_2}) =
H_n(\underline{\mu_3}) \oplus H_n(\underline{\mu_4}),
\]
where $\mu_4=\nu\tau_{i-1}$. And the last part of the lemma
follows.
\end{proof}

Similar to the proof of Proposition \ref{resolved-braids}, we can
induct on the weight of a resolved braid using Lemmas
\ref{induction}, \ref{hn-moves}, and prove:

\begin{proposition}\label{resolved-braids-n}
If $\mu=\tau_{i_1}\cdots\tau_{i_m}$ is a word in $\tau_{1}$,
$\cdots$, $\tau_{b-1}$, then
\[
-(n-1)b-m \leq g^{(n)}_{min}(\underline{\mu})\leq
g^{(n)}_{max}(\underline{\mu}) \leq (n-1)b+m.
\]
\end{proposition}
\begin{proof}
\end{proof}

\begin{proof}[Proof of inequality \eqref{main-inequality-n} in Theorem
\ref{main}] Let $B$ be a braid with $b$ strands $c_+$ positive
crossings, $c_-$ negative crossings, and writhe $w=c_+-c_-$. Note
that
\[
H_n(B)=H(H(C_n(B),d_{mf}),d_\chi),
\]
and
\[
H(C_n(B),d_{mf})=\bigoplus_{\Gamma}H_n(\Gamma)\{p_\Gamma\},
\]
where $\Gamma$ runs through all resolutions of $B$, and
$p_\Gamma\in \zed$. Let $\Gamma$ be any resolution of $B$. From
the construction of $C_n(B)$, it's clear that
\[
p_\Gamma=(n-1)w+c_{\Gamma,+}-c_{\Gamma,-},
\]
where $c_{\Gamma,+}$ (resp. $c_{\Gamma,-}$) is the number of
wide edges in $\Gamma$ from positive (resp. negative) crossings in
$B$. Let $g$ be the quantum degree of a non-vanishing
homogeneous element of $H_n(\Gamma)\{p_\Gamma\}$. Then
\[
g^{(n)}_{min}(\Gamma)+p_\Gamma \leq g \leq
g^{(n)}_{min}(\Gamma)+p_\Gamma.
\]
But, from Proposition \ref{resolved-braids-n}, we have
\[
-(n-1)b-c_{\Gamma,+}-c_{\Gamma,-}\leq g^{(n)}_{min}(\Gamma) \leq
g^{(n)}_{max}(\Gamma) \leq (n-1)b+c_{\Gamma,+}+c_{\Gamma,-}.
\]
So
\[
(n-1)(w-b)-2c_{\Gamma,-} \leq g \leq (n-1)(w+b)+2c_{\Gamma,+}.
\]
It's clear that $c_{\Gamma,+}\leq c_+$, and $c_{\Gamma,-}\leq
c_-$. This shows that, if $g$ is the quantum degree of any
non-vanishing homogeneous element of $H(C_n(B),d_{mf})$, then
\[
(n-1)(w-b)-2c_- \leq g \leq (n-1)(w+b)+2c_+.
\]
Thus,
\[
(n-1)(w-b)-2c_- \leq g^{(n)}_{min}(B) \leq g^{(n)}_{max}(B) \leq
(n-1)(w+b)+2c_+.
\]
\end{proof}

\section{Transversal Knot Invariants $\{\psi_n\}$}\label{inv}

In this section, we use the usual notation for braid group. Denote
by $\mathcal{B}_b$ the braid group of $b$ strands, and
$\sigma_1^{\pm1},\cdots,\sigma_{b-1}^{\pm1}$ the standard
generators of $\mathcal{B}_b$, i.e., $\sigma_i^{\pm1}$ permutes
the $i$-th and $(i+1)$-th strands with a $(\pm)$-crossing. Also,
denote by $\phi$ the empty word, which is the identity of
$\mathcal{B}_b$. Two braids represent the same link if and only if
one can be changed into the other by following operations:
\begin{itemize}
    \item Braid group relations generated by:
    $\sigma_i\sigma_i^{-1}=\sigma_i^{-1}\sigma_i=\phi$,
    $\sigma_i\sigma_j=\sigma_j\sigma_i$, when $|i-j|>1$, and
    $\sigma_i\sigma_{i+1}\sigma_i=\sigma_{i+1}\sigma_i\sigma_{i+1}$.
    \item Conjugations: $\mu\leftrightarrow\eta^{-1}\mu\eta$,
    where $\mu,~\eta\in \mathcal{B}_b$.
    \item Stabilizations and destabilizations:
    \[
    \left\{%
\begin{array}{ll}
    \text{positive}: & \mu~(\in \mathcal{B}_b)
    \leftrightarrow \mu\sigma_b~(\in \mathcal{B}_{b+1}) \\
    \text{negative}: & \mu~(\in \mathcal{B}_b)
    \leftrightarrow \mu\sigma_b^{-1}~(\in \mathcal{B}_{b+1}) \\
\end{array}%
\right.
\]
\end{itemize}

In \cite{Ben}, Bennequin proved that any transversal link is
transversally isotopic to a transversal braid. The following theorem
from \cite{OSh,Wr} describes when two transversal braids are
transversally isotopic.

\begin{theorem}[\cite{OSh,Wr}]\label{transversal-markov}
If two transversal braids are transversally isotopic, then one can be
changed into the other by a finite sequence of braid group
relations, conjugations and positive stabilizations and
destabilizations.
\end{theorem}

Let $\mu$ be a closed braid with $b$ strands, and $\Gamma$ the
resolution of $\mu$ so that each crossing of $\mu$ is $0$-resolved
as shown in Figure \ref{resolutions}. Then $\Gamma$ consists of
$b$ concentric circles, and
\[
H_n(\Gamma)= A_1 \otimes_{\mathbb{Q}} A_2 \otimes_{\mathbb{Q}}
\cdots \otimes_{\mathbb{Q}} A_b,
\]
where $A_i=[\mathbb{Q}[x_i]/(x_i^n)]\{-n+1\}$, and $\deg x_i =2$
as before. Let $w$ be the writhe of $\mu$. It's clear that the
standard isomorphism
\[
f: H_n(\Gamma) \xrightarrow{\cong} H_n(\Gamma)\{(n-1)w\} ~(\subset
H(C_n(\mu),d_{mf})),
\]
has quantum degree $(n-1)w$, and maps $H_n(\Gamma)$ into
cohomological degree $0$. Then the element
\[
\varphi_n(\mu)=f(x_1^{n-1} \otimes \cdots \otimes x_b^{n-1})
\]
has quantum degree $(n-1)(w+b)$ and cohomological degree $0$.

\begin{proposition}\label{phi-n}
$\varphi_n(\mu)$ is a cocycle in the chain complex
$(H(C_n(\mu),d_{mf}), d_\chi)$. And, up to multiplication by a
non-vanishing scalar, its cohomology class $[\varphi_n(\mu)]\in
H_n(\mu)=H(H(C_n(\mu),d_{mf}), d_\chi)$ is invariant under braid
group relations, conjugations and negative stabilizations and
destabilizations.
\end{proposition}
\begin{proof}
For a negative crossing $\varsigma$ of $\mu$, denote by
$\Gamma_\varsigma$ the resolution of $\mu$ obtained by
$1$-resolving $\varsigma$ and $0$-resolving all other crossings.
One can easily see that
\[
d_\chi(H_n(\Gamma)\{(n-1)w\}) \subset
\bigoplus_{\varsigma}H_n(\Gamma_\varsigma)\{(n-1)w-1\} \subset
H(C_n(\mu),d_{mf}),
\]
where $\varsigma$ runs through all negative crossings of $\mu$.
From the second part of Lemma \ref{hn-moves}, we know that the
maximal quantum degree of a non-vanishing homogeneous element of
$H_n(\Gamma_\varsigma)\{(n-1)w-1\}$ is $(n-1)(w+b)-2$. But
$d_\chi$ preserves the quantum grading. So
$d_\chi(\varphi_n(\mu))$ is a homogeneous element of degree
$(n-1)(w+b)$. This implies that $d_\chi(\varphi_n(\mu))=0$, and,
hence, $\varphi_n(\mu)$ is a cocycle.

Next, we prove the invariance of $[\varphi_n(\mu)]$ under the
local moves depicted in Figure \ref{braid-moves}, which implies
the second part of the proposition.

\begin{figure}[h]

\setlength{\unitlength}{1pt}

\begin{picture}(420,175)(-210,-130)


\put(-155,0){\vector(0,1){40}}

\put(-150,15){$\leftrightarrow$}

\put(-125,20){\vector(-1,2){10}}

\qbezier(-125,20)(-115,0)(-115,20)

\qbezier(-115,20)(-115,40)(-122.5,22.5)

\put(-135,0){\line(1,2){8.75}}

\put(-150,-30){\small{(i)}}

\put(-160,-15){$\mu_0$}

\put(-130,-15){$\mu_1$}


\put(-30,0){\vector(0,1){40}}

\put(-10,0){\vector(0,1){40}}

\put(-5,15){$\leftrightarrow$}

\put(10,30){\vector(0,1){10}}

\put(10,30){\line(4,-1){20}}

\put(30,15){\line(0,1){10}}

\put(10,10){\line(4,1){20}}

\put(10,0){\line(0,1){10}}

\put(30,30){\vector(0,1){10}}

\put(30,30){\line(-4,-1){7.5}}

\put(17.5,26.875){\line(-4,-1){7.5}}

\put(10,15){\line(0,1){10}}

\put(30,10){\line(-4,1){7.5}}

\put(17.3,13.125){\line(-4,1){7.5}}

\put(30,0){\line(0,1){10}}

\put(-5,-30){\small{(ii$_a$)}}

\put(-25,-15){$\mu_0$}

\put(15,-15){$\mu_1$}


\put(-30,-100){\vector(0,1){40}}

\put(-10,-100){\vector(0,1){40}}

\put(-5,-85){$\leftrightarrow$}

\put(10,-70){\vector(0,1){10}}

\put(10,-70){\line(4,-1){7.5}}

\put(30,-75){\line(-4,1){7.5}}

\put(30,-85){\line(0,1){10}}

\put(30,-85){\line(-4,-1){7.5}}

\put(10,-90){\line(4,1){7.5}}

\put(10,-100){\line(0,1){10}}

\put(30,-70){\vector(0,1){10}}

\put(30,-70){\line(-4,-1){20}}

\put(10,-85){\line(0,1){10}}

\put(30,-90){\line(-4,1){20}}

\put(30,-100){\line(0,1){10}}

\put(-5,-130){\small{(ii$_b$)}}

\put(-25,-115){$\mu_0$}

\put(15,-115){$\mu_1$}


\put(90,0){\vector(1,1){40}}

\put(100,0){\line(0,1){7.5}}

\put(100,12.5){\vector(0,1){27.5}}

\put(130,0){\line(-1,1){17.5}}

\put(107.5,22.5){\line(-1,1){5}}

\put(97.5,32.5){\vector(-1,1){7.5}}

\put(135,15){$\leftrightarrow$}

\put(150,0){\vector(1,1){40}}

\put(180,0){\line(0,1){27.5}}

\put(180,32.5){\vector(0,1){7.5}}

\put(190,0){\line(-1,1){7.5}}

\put(177.5,12.5){\line(-1,1){5}}

\put(167.5,22.5){\vector(-1,1){17.5}}

\put(135,-30){\small{(iii)}}

\put(105,-15){$\mu_0$}

\put(165,-15){$\mu_1$}

\end{picture}

\caption{Local moves}\label{braid-moves}

\end{figure}
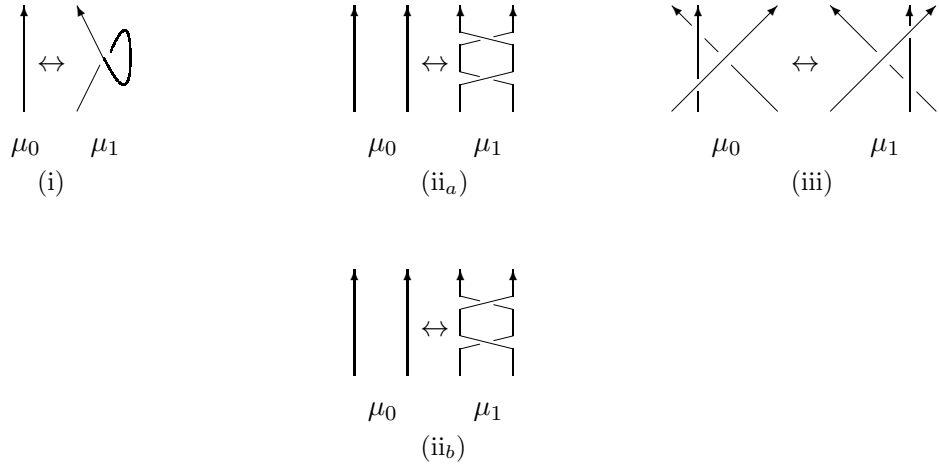

\begin{remark}
All together, there are $6$ braidlike Reidemeister type $3$ local moves to deal with. These correspond to the $6$ possible ordering of the three branches involved. But each of the other $5$ braidlike Reidemeister type $3$ local moves is related to move (iii) by braidlike Reidemeister type $2$ local moves. So we only need to prove the invariance under move (iii) here.
\end{remark}

\noindent\textbf{Invariance under move (i).} Consider the planar
diagrams in Figure \ref{braid-move-i}.
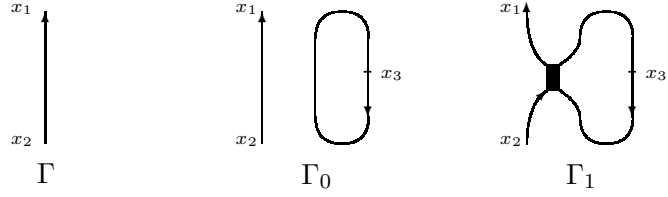
\begin{figure}[h]

\setlength{\unitlength}{1pt}

\begin{picture}(420,65)(-210,-25)

\linethickness{.5pt}


\put(-102,-10){\vector(0,1){50}}

\put(-115,-10){\tiny{$x_2$}}

\put(-115,40){\tiny{$x_1$}}

\put(-105,-25){$\Gamma$}


\put(-20,-10){\vector(0,1){50}}

\put(0,0){\line(0,1){30}}

\qbezier(0,30)(0,40)(10,40)

\qbezier(10,40)(20, 40)(20,30)

\put(20,30){\vector(0,-1){30}}

\qbezier(0,0)(0,-10)(10,-10)

\qbezier(10,-10)(20, -10)(20,0)

\put(-5,-25){$\Gamma_0$}

\put(-30,-10){\tiny{$x_2$}}

\put(-30,40){\tiny{$x_1$}}

\put(25,15){\tiny{$x_3$}}

\put(18.5,15){-}


\put(87.5,10){\vector(1,1){0}}

\qbezier(80,-10)(80,5)(87.5,10)

\qbezier(80,40)(80,25)(87,20)

\put(80,44){\vector(0,1){0}}

\qbezier(100,0)(100,5)(92.5,10)

\qbezier(100,30)(100,25)(92.5,20)

\qbezier(100,30)(100,40)(110,40)

\qbezier(110,40)(120, 40)(120,30)

\put(120,30){\vector(0,-1){30}}

\qbezier(100,0)(100,-10)(110,-10)

\qbezier(110,-10)(120, -10)(120,0)

\put(95,-25){$\Gamma_1$}

\put(70,-10){\tiny{$x_2$}}

\put(70,40){\tiny{$x_1$}}

\put(125,15){\tiny{$x_3$}}

\put(118.5,15){-}

\linethickness{5pt}

\put(90,10){\line(0,1){10}}

\end{picture}

\caption{Local diagrams related to move (i)}\label{braid-move-i}

\end{figure}

The local complex associated to $\mu_0$ is
\[
0 \rightarrow H_n(\Gamma) \rightarrow 0.
\]
And the local complex associated to $\mu_1$ is
\[
0 \rightarrow H_n(\Gamma_0)\{1-n\} \rightarrow H_n(\Gamma_1)\{-n\}
\rightarrow 0.
\]
Note that $H_n(\Gamma_0)= H_n(\Gamma)\otimes A_3$, where
$A_3=[\mathbb{Q}[x_3]/(x_3^n)]\{-n+1\}$, and $\deg x_3 =2$. There
is a homomorphism $A_3\xrightarrow{\varepsilon}\mathbb{Q}$ of
degree $1-n$ given by
\[
\varepsilon(x_3^j)=\left\{%
\begin{array}{ll}
    1, & \hbox{if $j=n-1$;} \\
    0, & \hbox{otherwise.} \\
\end{array}%
\right.
\]
This induces a quantum degree preserving homomorphism $\varepsilon:
H_n(\Gamma_0)\{1-n\}\rightarrow H_n(\Gamma)$, which, in turn,
induces a bi-degree preserving homomorphism of complexes
\[
\begin{array}{ccccccc}
  0 & \rightarrow & H_n(\Gamma_0)\{1-n\} & \rightarrow & H_n(\Gamma_1)\{-n\} & \rightarrow & 0 \\
   &  & \downarrow &  & \downarrow &  &  \\
  0 & \rightarrow & H_n(\Gamma) & \rightarrow & 0 &  &  \\
\end{array}.
\]
It's clear that this chain map maps $\varphi_n(\mu_1)$ to
$\varphi_n(\mu_0)$. And, from Section 8 of \cite{KR1}, we know that
this chain map induces (a non-zero multiple of) the standard
isomorphism from $H_n(\mu_1)$ to $H_n(\mu_0)$. Thus,
$[\varphi_n(\mu)]$ is invariant under local move (i) up to
multiplication by non-vanishing scalar.

\noindent\textbf{Invariance under moves (ii$_a$) and (ii$_b$).} The
proofs for these two moves are almost identical. So we only prove
the invariance under move (ii$_a$).

Consider the planar diagrams in Figure \ref{braid-move-ii}.
\begin{figure}[h]

\setlength{\unitlength}{1pt}

\begin{picture}(420,160)(-210,-10)

\linethickness{.5pt}


\put(-70,80){\vector(0,1){60}}

\put(-50,80){\vector(0,1){60}}

\put(-100,105){$\Gamma_{00}$}

\put(-80,140){\tiny{$x_1$}}

\put(-80,80){\tiny{$x_4$}}

\put(-47,140){\tiny{$x_2$}}

\put(-47,80){\tiny{$x_3$}}


\linethickness{.5pt}

\put(62.5,120){\vector(-3,2){0}}

\qbezier(70,80)(70,120)(62.5,120)

\put(57.5,120){\vector(3,2){0}}

\qbezier(50,80)(50,120)(57.5,120)

\qbezier(50,140)(50,135)(57.5,130)

\qbezier(70,140)(70,135)(62.5,130)

\put(70,140){\vector(1,4){0}}

\put(50,140){\vector(-1,4){0}}

\put(20,105){$\Gamma_{10}$}

\put(40,140){\tiny{$x_1$}}

\put(40,80){\tiny{$x_4$}}

\put(73,140){\tiny{$x_2$}}

\put(73,80){\tiny{$x_3$}}

\linethickness{5pt}

\put(60,120){\line(0,1){10}}


\linethickness{.5pt}

\put(62.5,10){\vector(-1,1){0}}

\qbezier(70,0)(70,5)(62.5,10)

\qbezier(70,30)(70,25)(62.5,20)

\put(70,30){\vector(1,4){0}}

\put(50,30){\vector(-1,4){0}}

\put(57.5,10){\vector(1,1){0}}

\qbezier(50,0)(50,5)(57.5,10)

\qbezier(50,30)(50,25)(57.5,20)

\qbezier(70,30)(70,35)(62.5,40)

\qbezier(70,60)(70,55)(62.5,50)

\put(70,60){\vector(1,4){0}}

\put(50,60){\vector(-1,4){0}}

\qbezier(50,30)(50,35)(57.5,40)

\qbezier(50,60)(50,55)(57.5,50)

\put(20,25){$\Gamma_{11}$}

\put(40,60){\tiny{$x_1$}}

\put(40,0){\tiny{$x_4$}}

\put(73,60){\tiny{$x_2$}}

\put(73,0){\tiny{$x_3$}}

\put(40,30){\tiny{$x_6$}}

\put(73,30){\tiny{$x_5$}}

\put(48.5,30){-}

\put(68.5,30){-}

\linethickness{5pt}

\put(60,10){\line(0,1){10}}

\put(60,40){\line(0,1){10}}


\linethickness{.5pt}

\put(-62.5,10){\vector(1,1){0}}

\qbezier(-70,0)(-70,5)(-62.5,10)

\put(-57.5,10){\vector(-1,1){0}}

\qbezier(-50,0)(-50,5)(-57.5,10)

\qbezier(-70,60)(-70,20)(-62.5,20)

\put(-70,60){\vector(0,0){0}}

\put(-50,60){\vector(0,1){0}}

\qbezier(-50,60)(-50,20)(-57.5,20)

\put(-100,25){$\Gamma_{01}$}

\put(-80,60){\tiny{$x_1$}}

\put(-80,0){\tiny{$x_4$}}

\put(-47,60){\tiny{$x_2$}}

\put(-47,0){\tiny{$x_3$}}

\linethickness{5pt}

\put(-60,10){\line(0,1){10}}


\linethickness{.5pt}

\put(-172.5,50){\vector(1,1){0}}

\qbezier(-180,40)(-180,45)(-172.5,50)

\put(-167.5,50){\vector(-1,1){0}}

\qbezier(-160,40)(-160,45)(-167.5,50)

\qbezier(-180,100)(-180,95)(-172.5,90)

\put(-180,100){\vector(-1,4){0}}

\put(-160,100){\vector(1,4){0}}

\qbezier(-160,100)(-160,95)(-167.5,90)

\put(-173,30){$\Gamma$}

\put(-190,100){\tiny{$x_1$}}

\put(-190,40){\tiny{$x_4$}}

\put(-157,100){\tiny{$x_2$}}

\put(-157,40){\tiny{$x_3$}}

\linethickness{5pt}

\put(-170,50){\line(0,1){40}}

\end{picture}

\caption{Local diagrams related to move
(ii$_a$)}\label{braid-move-ii}
\end{figure}
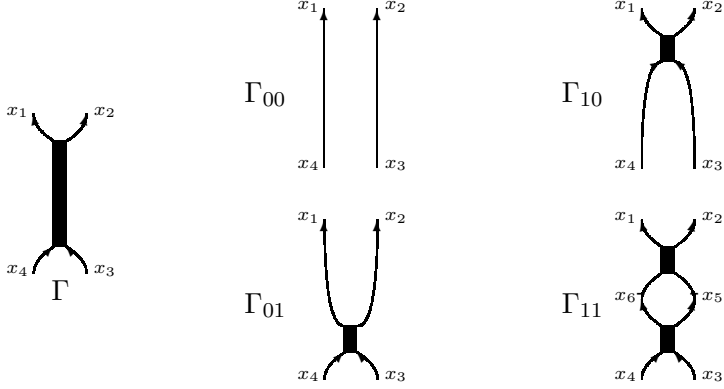

The local complex associated to $\mu_0$ is
\[
0 \rightarrow H_n(\Gamma_{00}) \rightarrow 0.
\]
And the local complex associated to $\mu_1$ is
\[
0 \rightarrow H_n(\Gamma_{01})\{1\} \rightarrow H_n(\Gamma_{00})
\oplus H_n(\Gamma_{11}) \rightarrow H_n(\Gamma_{10})\{-1\}
\rightarrow 0.
\]
Let
\[
\pi_{11}: H_n(\Gamma_{00}) \oplus H_n(\Gamma_{11}) \rightarrow
H_n(\Gamma_{11})
\]
be the natural projection. Consider the quantum degree preserving
homomorphisms
\[
f_1 = \pi_{11}\circ d_\chi: H_n(\Gamma_{01})\{1\} \rightarrow
H_n(\Gamma_{11})
\]
and
\[
f_2 = d_\chi|_{H_n(\Gamma_{11})}: H_n(\Gamma_{11}) \rightarrow
H_n(\Gamma_{10})\{-1\}.
\]
Note that
\[
H_n(\Gamma_{11}) \cong H_n(\Gamma)\{1\} \oplus H_n(\Gamma)\{-1\}.
\]
Let
\[
f_{10}:H_n(\Gamma_{01})\{1\} \rightarrow H_n(\Gamma)\{1\} ~\text{
and }~ f_{21}:H_n(\Gamma)\{-1\} \rightarrow H_n(\Gamma_{10})\{-1\}
\]
be the obvious maps induced by $f_1$, $f_2$ and the above
isomorphism. Note that $\Gamma_{01}$ and $\Gamma_{10}$ are isotopic
to $\Gamma$. In Section 8 of \cite{KR1}, M. Khovanov and L. Rozansky
showed that $f_{10}$ and $f_{21}$ are non-zero multiples of the
natural isomorphisms
\[
H_n(\Gamma_{01})\{1\} \xrightarrow{\cong} H_n(\Gamma)\{1\} ~\text{
and }~ H_n(\Gamma)\{-1\} \xrightarrow{\cong} H_n(\Gamma_{10})\{-1\}.
\]
Define submodule $A \subset H_n(\Gamma_{00}) \oplus
H_n(\Gamma)\{-1\} \subset H_n(\Gamma_{00})
\oplus H_n(\Gamma_{11})$ by
\[
A=\{(a,-f_{21}^{-1} \circ d_\chi(a))|~a\in H_n(\Gamma_{00})\},
\]
and submodule $B\subset H_n(\Gamma_{00})
\oplus H_n(\Gamma_{11})$ by
\[
B= d_\chi (H_n(\Gamma_{01})\{1\}).
\]
Then
\[
H_n(\Gamma_{00}) \oplus H_n(\Gamma_{11}) = A \oplus B \oplus
H_n(\Gamma)\{-1\}.
\]
It's easy to see that $f_{21}=d_\chi|_{H_n(\Gamma)\{-1\}}$. Thus,
the local complex associated to $\mu_1$ splits into the direct sum
of the following chain complexes
\[
\begin{array}{ccccccccc}
   &  & 0 & \rightarrow & A & \rightarrow & 0, & &  \\
  0 & \rightarrow & H_n(\Gamma_{01})\{1\}\} & \xrightarrow{\cong} & B & \rightarrow & 0, &  & \\
   &  & 0 & \rightarrow & H_n(\Gamma)\{-1\} & \xrightarrow{\cong} & H_n(\Gamma_{10})\{-1\} & \rightarrow &
0.
\end{array}
\]
The second and third chain complexes give rise to trivial
cohomologies. And the standard isomorphism from $H_n(\mu_0)$ to $H_n(\mu_1)$
is induced by (a non-zero scalar multiple of) the isomorphism of
chain complexes
\[
F: [0 \rightarrow H_n(\Gamma_{00}) \rightarrow 0]
\xrightarrow{\cong} [0 \rightarrow A \rightarrow 0]
\]
given by
\[
F(a)=(a,-f_{21}^{-1} \circ d_\chi(a)).
\]
Since $d_\chi(\varphi_n(\mu_0))=0$, it's easy to see that
$F(\varphi_n(\mu_0))=\varphi_n(\mu_1)$. Thus, $[\varphi_n(\mu)]$ is
invariant under local move (ii$_a$) up to multiplication by
non-vanishing scalar.

\noindent\textbf{Invariance under moves (iii).}
The proof in \cite{KR1} of the invariance of $H_n$ under move
(iii) is much more complex than that under moves
(i) and (ii). But, if we first use D. Bar-Natan's algebraic trick in
\cite{Bar} to reformulate the description of the isomorphism of
cohomologies, then it is quite easy to establish the invariance of
$[\varphi_n(\mu)]$ under this move.

\begin{figure}[h]

\setlength{\unitlength}{1pt}

\begin{picture}(420,160)(-210,-80)


\qbezier(-90,20)(-80,30)(-80,80)

\put(-80,80){\vector(0,1){0}}

\qbezier(-75,20)(-75,40)(-35,80)

\put(-35,80){\vector(1,1){0}}

\put(-30,20){\line(-1,1){29}}

\put(-63,53){\line(-1,1){14.5}}

\put(-82.5,72.5){\vector(-1,1){7.5}}

\put(-110,50){$\Gamma_{00}$}


\qbezier(80,20)(80,70)(90,80)

\put(90,80){\vector(3,4){0}}

\qbezier(35,20)(75,60)(75,80)

\put(75,80){\vector(0,1){0}}

\put(90,20){\line(-1,1){7.5}}

\put(77.5,32.5){\line(-1,1){14.5}}

\put(58,52){\vector(-1,1){28}}

\put(10,50){$\Gamma_{10}$}


\put(-80,-60){\vector(0,1){40}}

\put(-75,-80){\vector(0,1){10}}

\put(-90,-80){\vector(1,1){10}}

\put(-75,-60){\vector(1,1){40}}

\put(-30,-80){\line(-1,1){30}}

\put(-65,-45){\line(-1,1){12.5}}

\put(-82.5,-27.5){\vector(-1,1){7.5}}

\put(-110,-50){$\Gamma_{01}$}


\put(80,-80){\vector(0,1){40}}

\put(80,-30){\vector(1,1){10}}

\put(75,-30){\vector(0,1){10}}

\put(35,-80){\vector(1,1){40}}

\put(90,-80){\line(-1,1){7.5}}

\put(77.5,-67.5){\line(-1,1){12.5}}

\put(60,-50){\vector(-1,1){30}}

\put(10,-50){$\Gamma_{11}$}

\linethickness{5pt}

\put(-77.5,-70){\line(0,1){10}}

\put(77.5,-40){\line(0,1){10}}

\end{picture}

\caption{Local diagrams related to move
(iii)}\label{braid-move-iii-a}

\end{figure}
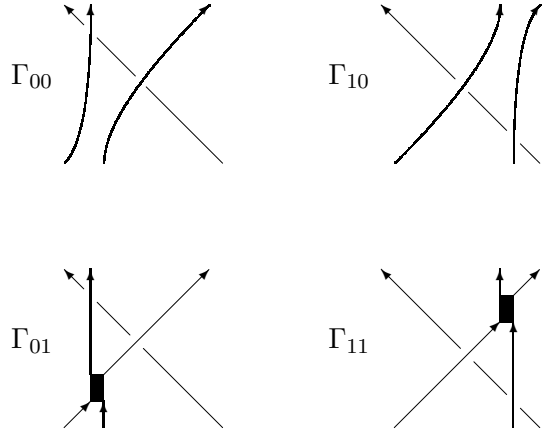

Consider the partially resolved planar diagrams in Figure \ref{braid-move-iii-a}. We
extend the definition of the (double) chain complex
$(C_n,d_{mf},d_\chi)$ to these diagrams. Then, for $i=0,1$, the chain complex
$(H(C_n(\mu_i),d_{mf}),d_\chi)$ is the mapping cone of the chain map
\[
\chi_{i1}: (H(C_n(\Gamma_{i1}),d_{mf}),d_\chi) \rightarrow (H(C_n(\Gamma_{i0}),d_{mf}),d_\chi),
\]
where $\chi_{i1}$ is induced by the $\chi^{(n)}_1$-map
associated to the wide edge in $\Gamma_{i1}$. Note that $\Gamma_{00}$ and $\Gamma_{10}$ are the same diagram, which we denote by $\Gamma$. Although formulated differently in \cite{KR1}, M. Khovanov and L. Rozansky actually constructed a chain complex $\mathfrak{C}$ which is a strong deformation retract of both $(H(C_n(\Gamma_{01}),d_{mf}),d_\chi)$ and $(H(C_n(\Gamma_{11}),d_{mf}),d_\chi)$ with corresponding inclusion maps 
\[
F_i: \mathfrak{C} \rightarrow (H(C_n(\Gamma_{i1}),d_{mf}),d_\chi), ~i=0,1,
\] 
so that $\chi_{01}\circ F_0 = \chi_{11} \circ F_1$ as homomorphisms from $\mathfrak{C}$ to $(H(C_n(\Gamma),d_{mf}),d_\chi)$.
Thus, the following squares commute:
\[
\begin{array}{ccccc}
  (H(C_n(\Gamma_{00}),d_{mf}),d_\chi) & \xrightarrow{\id} & (H(C_n(\Gamma),d_{mf}),d_\chi) & \xleftarrow{\id} & (H(C_n(\Gamma_{10}),d_{mf}),d_\chi) \\
  \tiny{\chi_{01}}\uparrow & &\tiny{\chi}\uparrow & & \tiny{\chi_{11}}\uparrow \\
  (H(C_n(\Gamma_{01}),d_{mf}),d_\chi) &\xleftarrow{F_0} & \mathfrak{C} & \xrightarrow{F_1} & (H(C_n(\Gamma_{11}),d_{mf}),d_\chi)
\end{array},
\]
where $\id$ is the natural identity homomorphism, and $\chi=\chi_{01}\circ F_0 = \chi_{11} \circ F_1$. 

By Lemma 4.5 of \cite{Bar}, for $i=0,1$, the maps $\id$ and $F_i$ induce an isomorphism $\mathcal{F}_i$ from $H_n(\mu_i)$ to the cohomology of the mapping cones of $\chi: \mathfrak{C} \rightarrow (H(C_n(\Gamma),d_{mf}),d_\chi)$. Then the isomorphism
\[
\mathcal{F}=\mathcal{F}_1^{-1} \circ \mathcal{F}_0:H_n(\mu_0)\rightarrow H_n(\mu_1)
\]
is (a multiple of) the standard isomorphism. It's clear that
\[
\varphi(\mu_0)\in H(C_n(\Gamma_{00}),d_{mf}) \subset
H(C_n(\mu_0),d_{mf}),
\]
\[
\varphi(\mu_1)\in H(C_n(\Gamma_{10}),d_{mf}) \subset
H(C_n(\mu_1),d_{mf}),
\]
and
\[
\id(\varphi(\mu_0))=\varphi(\mu_1),
\]
where $\id:H(C_n(\Gamma_{00}),d_{mf})\rightarrow H(C_n(\Gamma_{10}),d_{mf})$ is the natural identity homomorphism. This implies that
\[
\mathcal{F}([\varphi(\mu_0)]) = [\varphi(\mu_1)],
\]
which shows that $[\varphi_n(\mu)]$ is invariant under local move
(iii) up to multiplication by non-vanishing scalar.
\end{proof}

Let $L$ be a transversal link in the standard contact $S^3$. We
transversally isotope it into a transversal braid.

\begin{definition}
$\psi_n(L)=[\varphi_n(\overline{L})]\in H_n(\overline{L})$.
\end{definition}

\begin{proof}[Proof of Theorem \ref{invariants}]
From Theorem \ref{transversal-markov} and Proposition \ref{phi-n},
we know that, up to multiplication by non-zero scalar, $\psi_n(L)$
is independent of the transversal braid representation, and is
therefore invariant under transversal isotopy. From the
construction of $\psi_n(L)=[\varphi_n(\overline{L})]$, it's clear
that it has cohomological degree $0$ and quantum degree
$(n-1)(-w+b)=-(n-1)sl(L)$, where $w$ is the writhe of $L$, and $b$
is the number of strands of $L$. Compare our construction with
that by Plamenevskaya in \cite{Pl4}. One easily sees that
$\psi_2(L)$ is identified with $\psi(L)$ under the isomorphism
$H_2(\overline{L})\cong \mathcal{H}(L)\otimes_{\zed}\mathbb{Q}$.
\end{proof}

Next we generalize Plamenevskaya's result about quasi-positive
transversal braids, which completes the proof of Corollary
\ref{top-sl}.

\begin{proposition}[Compare Theorem 4 of
\cite{Pl4}]\label{quasi-positive} Let $\mu\in \mathcal{B}_b$ be a
closed braid, and $\mu'=\mu\sigma_i^{-1}$, where $1\leq i \leq
b-1$. Then there is a homomorphism $f:H_n(\mu')\rightarrow
H_n(\mu)$ of quantum degree $n-1$ such that
$f([\varphi_n(\mu')])=[\varphi_n(\mu)]$. Specially, this implies
that the $\psi_n$ invariant of a quasi-positive transversal braid
is non-vanishing.
\end{proposition}

\begin{proof}

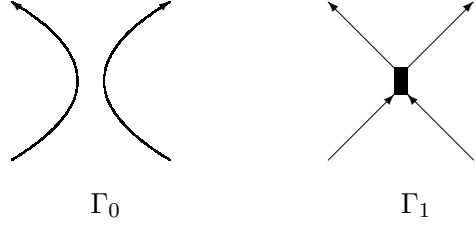
\begin{figure}[h]

\setlength{\unitlength}{1pt}

\begin{picture}(420,80)(-210,0)


\qbezier(-90,20)(-40,50)(-90,80)

\put(-90,80){\vector(-4,3){0}}

\qbezier(-30,20)(-80,50)(-30,80)

\put(-30,80){\vector(4,3){0}}

\put(-60,0){$\Gamma_0$}


\put(85,20){\vector(-1,1){25}}

\put(30,20){\vector(1,1){25}}

\put(60,55){\vector(1,1){25}}

\put(55,55){\vector(-1,1){25}}

\linethickness{5pt}

\put(57.5,45){\line(0,1){10}}

\put(57.5,0){$\Gamma_1$}

\end{picture}

\caption{Diagrams $\Gamma_0$ and $\Gamma_1$}\label{remove+}

\end{figure}

The local chain complex associated to $\sigma_i^{-1}$ in $(H_n(\mu',d_{mf}),d_\chi)$
is
\[
0 \rightarrow H_n(\Gamma_0)\{1-n\} \rightarrow H_n(\Gamma_1)\{-n\}
\rightarrow 0.
\]
And the corresponding local complex in $(H_n(\mu,d_{mf}),d_\chi)$
is
\[
0 \rightarrow H_n(\Gamma_0) \rightarrow 0.
\]
The identity map of $H_n(\Gamma_0)$ induces a local chain map of
quantum degree $n-1$:
\[
\begin{array}{ccccccc}
  0 & \rightarrow & H_n(\Gamma_0)\{1-n\} & \rightarrow & H_n(\Gamma_1)\{-n\} & \rightarrow & 0 \\
   &  & \downarrow &  & \downarrow &  &  \\
  0 & \rightarrow & H_n(\Gamma_0) & \rightarrow & 0 &  &  \\
\end{array}.
\]
Let $f:H_n(\mu')\rightarrow H_n(\mu)$ be the homomorphism induced
by the above local chain map. It's clear that
$f([\varphi_n(\mu')])=[\varphi_n(\mu)]$.

Now let $L$ be a quasi-positive transversal braid of $b$ strands,
i.e., $L$ is represented by a braid word of the form
\[
\nu=\mu_1\sigma_{i_1}\mu_1^{-1}\cdots\mu_k\sigma_{i_k}\mu_k^{-1}.
\]
Then $\overline{L}$ is represented by the word
\[
\overline{\nu} =
\overline{\mu}_1\sigma_{i_1}^{-1}\overline{\mu}_1^{-1} \cdots
\overline{\mu}_k\sigma_{i_k}^{-1}\overline{\mu}_k^{-1},
\]
where $\overline{\mu}_j$ is the mirror image of $\mu_j$. From the
first half of the proposition and the invariance of $[\varphi_n]$,
we know that there is a homomorphism $F:H_n(\overline{L})
\rightarrow H_n(\phi)$ such that $F([\varphi_n(\overline{\nu})]) =
[\varphi_n(\phi)]\neq 0$, where $\phi$ is the empty word in $b$
strands. Thus, $\psi_n(L)=[\varphi_n(\overline{\nu})]\neq0$.
\end{proof}

Some other properties of Plamenevskaya's $\psi$ invariant
generalize to $\psi_n$ too. For example:

\begin{proposition}[Compare Proposition 3 of
\cite{Pl4}]\label{neg-stabilization} Let $\mu\in \mathcal{B}_b$ be
a braid with $b$ strands. If there is an $i\in\{1,\cdots,b-1\}$
such that $\sigma_i$ occurs in the word $\mu$, but $\sigma_i^{-1}$
does not, then $[\varphi_n(\mu)]=0$, $\forall ~n\geq 2$.
Specially, if a transversal link $L$ is a transversal
stabilization of another transversal link, then $\psi_n(L)=0$,
$\forall ~n\geq 2$.
\end{proposition}

\begin{proof}
Fix a crossing $\sigma_i$ in $\mu$. Let $\Gamma_{-1}$ be the
resolution of $\mu$ obtained by $1$-resolving $\sigma_i$ and
$0$-resolving all other crossings, and $\Gamma_0$ the resolution
of $\mu$ obtained by $0$-resolving all the crossings. For a
negative crossing $\varsigma$ of $\mu$, let $\Gamma_{\varsigma}$
be the resolution of $\mu$ obtained by $1$-resolving $\sigma_i$,
$\varsigma$ and $0$-resolving all other crossings. Since $\mu$
does not contain $\sigma_i^{-1}$, $\Gamma_{\varsigma}$ is a
resolved closed braid of the form $\underline{\tau_i\tau_j}$, where
$j\neq i$. Then, in the chain complex
$(H(C_n(\mu),d_{mf}),d_\chi)$, we have
\[
d_\chi(H_n(\Gamma_{-1})\{(n-1)w+1\}) \subset
H_n(\Gamma_0)\{(n-1)w\} \bigoplus
(\bigoplus_{\varsigma}H_n(\Gamma_\varsigma)\{(n-1)w\}),
\]
where $w$ is the writhe of $\mu$, and $\varsigma$ runs through all
negative crossings of $\mu$. Let
\[
\pi_0: H_n(\Gamma_0)\{(n-1)w\} \bigoplus
(\bigoplus_{\varsigma}H_n(\Gamma_\varsigma)\{(n-1)w\}) \rightarrow
H_n(\Gamma_0)\{(n-1)w\}
\]
be the natural projection. Then it's clear that
\[
\pi_0 \circ d_\chi = \chi_1: H_n(\Gamma_{-1})\{(n-1)w+1\}
\rightarrow H_n(\Gamma_0)\{(n-1)w\},
\]
where $\chi_1$ is induced by the $\chi^{(n)}_1$-map related to the
crossing $\sigma_i$. From Proposition 29 of \cite{KR1}, we know
that there is a non-zero homogeneous element $\eta\in
H_n(\Gamma_{-1})\{(n-1)w+1\}$ of quantum degree $(n-1)(w+b)$ such
that $\chi_1(\eta)= \varphi_n(\mu) \in H_n(\Gamma_0)\{(n-1)w\}$.
By (1) and (2) of Lemma \ref{hn-moves}, we have that
$g^{(n)}_{max}(\Gamma_\varsigma)=(n-1)b-2$ for any negative
crossing $\varsigma$ of $\mu$. So the maximal quantum degree of a
non-zero homogeneous element of $H_n(\Gamma_\varsigma)\{(n-1)w\}$
is $(n-1)(w+b)-2$. Since $d_\chi$ preserves the quantum degree,
the projection of $d_\chi(\eta)$ onto
$H_n(\Gamma_\varsigma)\{(n-1)w\}$ is zero for any negative
crossing $\varsigma$ of $\mu$. This implies that $d_\chi(\eta)=
\varphi_n(\mu)$. Thus, $\varphi_n(\mu)$ is a coboundary, and
$[\varphi_n(\mu)]=0$.

If $L$ is a transversal stabilization of another transversal link,
then there is a $b\geq2$, such that $L$ is transversally isotopic
to a transversal braid given by a braid word $\mu \in
\mathcal{B}_b$, where $\mu$ contains one $\sigma_{b-1}^{-1}$ and
no $\sigma_{b-1}$. Therefore,
$\psi_n(L)=[\varphi_n(\overline{\mu})]=0$.
\end{proof}

\begin{acknowledgments}
The author's interest in finding new bounds for the self-linking
number was inspired by Matt Hedden's work on Legendrian knots and
the knot Floer homology. The basic idea of Theorem \ref{main} was
conceived during the Holomorphic Curves Workshop held at IAS in
May and June of 2005. The author would like to thank the
organizers for their hospitality, and Matt Hedden for interesting
conversations. The author would also like to thank the referee for many helpful comments on an earlier draft of this paper.
\end{acknowledgments}

\end{document}